\documentclass[12pt]{geophysics}

\usepackage{natbib}
\usepackage{graphicx}
\usepackage{color}
\usepackage{amsmath}
\usepackage{amsthm}
\usepackage{amssymb}
\usepackage{amsbsy}
\usepackage{setspace}
\usepackage{bm}
\usepackage{float}
\usepackage{hyperref}
\usepackage{algorithm}
\usepackage{algpseudocode}

\usepackage{placeins}
\newcommand{\bx}{\mathbf{x}}

\setfigdir{Fig}

\begin{document}
\title{A numerical investigation of Matched Source Waveform Inversion applied to acoustic transmission data}
\author{William W. Symes$^1$,  Huiyi Chen$^2$ and Susan E. Minkoff$^3$ \\
  $^1$ Department of Computational and Applied Mathematics, Rice University, Houston TX 77251 USA \\
  $^2$ Department of Mathematical Sciences, FO 35, 
  University of Texas at Dallas, Richardson, TX 75080 USA \\
  $^3$ Brookhaven National Laboratory, PO Box 5000, Upton, NY 11973-5000 USA
  }
\lefthead{Symes, Chen, Minkoff}

\righthead{Numerical Investigation of MSWI}

\maketitle
\parskip 12pt

\begin{abstract}
Iterative inversion of seismic, ultrasonic, and other wave data by local gradient-based optimization of mean-square data prediction error (Full Waveform Inversion or FWI) can fail to converge to useful model estimates if started from an initial model predicting wave arrival times in error by more than half a wavelength (a phenomenon known as cycle skipping). Matched Source Waveform Inversion (MSWI) extends the wave propagation model by a filter that shifts predicted waves to fit observed data. The MSWI objective adds a penalty for deviation of this filter from the identity to the mean-square data misfit . The extension allows the inversion to make large model adjustments while maintaining data fit and so reduces the chances of local optimization iterates stagnating at non-informative model estimates. Theory suggests that MSWI applied to acoustic transmission data with single-arrival wavefronts may produce an estimate of refractive index similar to the result of travel time inversion, but without requiring explicit identification of travel times. Numerical experiments conform to this expectation, in that MSWI applied to single arrival transmission data gives reasonable model estimates in cases where FWI fails. This MSWI model can then be used to jumpstart FWI for further refinement of the model. The addition of moderate amounts of noise (30\%) does not negatively impact MSWI's ability to converge. However, MSWI applied to data with multiple arrivals is no longer theoretically equivalent to travel-time tomography and exhibits the same tendency to cycle-skip as does FWI. 
\end{abstract}
\setlength{\parindent}{0cm}

\section{Introduction}
Matched Source Waveform Inversion (MSWI) \cite[]{HuangSymes2015SEG,HuangSymes:Geo17} is a variant of Full Waveform
Inversion (FWI) \cite[]{VirieuxOperto:09}, or least-squares minimization of the misfit between observed and predicted propagating wave data. The large size of the computational domain for typical FWI problems means that efficient local (Newton-like) optimization methods offer the only practical solution approaches. However, such algorithms have a tendency to stagnate at suboptimal model estimates, which fail to predict wave arrivals at correct times. Wave data is oscillatory, and the failure to predict correct arrival times manifests as matching of incorrect parts of the oscillatory wave data (``cycle-skipping'') \cite[]{VirieuxOperto:09}. 

This paper reviews the theoretical basis of MSWI, 
describes a practical computational framework, and 
numerically illustrates the theory by applying MSWI to solve several synthetic acoustic inversion problems for cases in which least-squares FWI fails. The theory suggests that for data generated from smooth acoustic models with single-arrival wavefronts, approximate minimization of the MSWI objective via a modest number of local optimization iterations can produce an intermediate model that then allows iterative FWI to successfully approximate a global minimizer.  
We present several numerical experiments supporting this contention. We show other examples in which (1) MSWI produces a smooth initial estimate from which FWI succeeds, despite non-smooth features in the target model; and (2) MSWI produces an initial estimate from which FWI succeeds, despite the presence of significant coherent data noise; and (3) MSWI fails to produce a suitable initial model for FWI due to the presence of multiple arrival times between source and receiver points. In all of these examples, FWI fails due to cycle skipping when
started from the same initial model as MSWI. 

MSWI loosens the bond between predicted and observed data by interposing a filter, adapted to map the two datasets trace-by-trace.
The MSWI objective function is a linear combination of the mean-square error between the observed data and the filtered predicted data and a quadratic penalty that vanishes when the filter is an impulse. 
This construction effectively relaxes the FWI data-fitting problem and permits local optimization methods to make large changes to the model and predicted data while maintaining fit of filtered predicted and observed data. 

The theoretical justification for MSWI lies in its link to travel time tomography \cite[]{Song:94c,Symes:94c,HuangSymes2015SEG,HuangSymes:Geo17,Symes:24a}. 
For smooth acoustic models producing simple wavefronts (single arrival times between sources and receivers), the MSWI objective function is proportional to a version of the travel time tomography objective, up to an error that decreases with increasing data frequency in time. Travel time tomography (minimization of the mean-square error between predicted and observed travel times) is commonly used to construct initial models for FWI  \cite[]{Bordingetal:87,SirguePratt:04,VirieuxOperto:09}. It is widely believed to be free of stationary points other than global minimizers, although theoretical support for this contention seems to be lacking. Thus we expect similar behaviour of the MSWI objective function, provided that the observed data is produced from a smooth model with single arrivals. Since matching of arrival times (to within a half-wavelength or so) is the conventional criterion for successful initiation of FWI, we expect that an approximate minimizer of the MSWI objective function may be a good initial guess for iterative FWI. The theory also suggests that the MSWI-tomography link may fail for multi-arrival data.

MSWI in the form described here is mathematically
equivalent to an {\em extended} formulation of the inverse problem, in
which acoustic point sources are allowed to depend on the positions of receiving
sensors, a non-physical expansion of the simulation domain. A number
of other variants of FWI have been based on essentially the same
extension of acoustic modeling
\cite[]{Song:94c,Symes:94c,Plessix:00,LuoSava:11,LiAlkhalifah:21}. In
particular this extension is also the basis of Adaptive Waveform Inversion (AWI) \cite[]{Warner:16,
  GuaschWarnerRavaut:GEO19,Warneretal:SEG21, Guaschetal:NPJDM20}. MSWI
is closely related to AWI but not identical. AWI includes a
normalization of the adaptive filter, which makes the AWI objective
function an even closer approximation to travel time mean square error
than is MSWI, at least for transmitted wave data with a single
arriving wavefront \cite[]{Symes:24a}.

\cite{HuangSymes:Geo17} also present a numerical exploration of MSWI, illustrating some of the same major points made in this paper. However, their work differs in two key respects from that reported here. First, \cite{HuangSymes:Geo17} use frequency-domain simulation of acoustic wave motion, whereas we use time-domain methods. Second, and more important, we use a penalty function formulation of MSWI, whereas \cite{HuangSymes:Geo17} use a zero-weight limit of this penalty function as the objective (the relation between the two approaches is explained in the Theory section below). Our previous work on a single-trace analogue of MSWI \cite[]{SymesMinkoffChen:IP22,ChenSymesMinkoff:IMAGE22} explored the potential benefits of the penalty approach, including the possibility of data noise estimation.


The next section describes versions of FWI and MSWI based on
acoustic wave propagation and reviews the theoretical connection between MSWI and
travel time inversion. The third section details the computational methods used to simulate acoustic data, apply local
optimization to the objective functions, implement variable projection reduction of MSWI, and control smoothness of the estimated model via
weighted norms in the domain of the simulation operator. In the main section of the paper we provide numerical examples that illustrate the theory.
The final Discussion Section describes several additional topics for future work including penalty weight selection, inclusion of more complex physics, and other extension methods that may overcome FWI's stagnation tendencies even with complex data kinematics.

\section{Theory}
For all the experiments described in this paper, we assume acoustic wave propagation with
isotropic point sources and receivers.
The pressure and velocity fields $p({\bf x},t;{\bf x}_s)$, ${\bf v}({\bf x},t;{\bf x}_s)$ for the source location ${\bf x}_s$ depend on the bulk modulus $\kappa({\bf x})$, buoyancy $\beta$ (reciprocal of the density $\rho$), and source wavelet $w(t)$ through the acoustic system
\begin{eqnarray}
  \label{eqn:awe}
 \frac{\partial p}{\partial t} & = &- \kappa \nabla \cdot {\bf v} +
\left(\int^t w\right) \delta({\bf x}-{\bf x}_s); \nonumber \\
\frac{\partial {\bf v}}{\partial t} & = & - \beta \nabla p; \\ 
p, {\bf v} & = & 0 \mbox{ for }  t \ll 0.
\end{eqnarray}
The model vectors, $m=(\kappa,\rho),$ and the source wavelet, $w,$ make up the domain of the forward
map or {\em modeling operator} $F[m] = \{p({\bf x}_r,t;{\bf
  x}_s)\}$, for specified source and receiver positions ${\bf x}_s, {\bf x}_r$ and recording
time interval.

In the discussion that follows, the buoyancy, $\beta,$ and the source wavelet, $w,$ will be regarded
as fixed, known parameters. Thus $m$ is effectively the bulk modulus field $\kappa$. In practice, neither buoyancy nor source characteristics may be known precisely and instead should be added to the parameters to be determined by inversion. 

\subsection{FWI and the Weighted Gradient}
In this setting FWI means: given buoyancy, 
source wavelet,
and data traces, $d(\bx_r,\cdot;\bx_s)$, find a model, $m,$ 
so that $F[m] \approx d$. The
simplest version of FWI 
involves finding a model
that minimizes the mean square error
\begin{equation}
  \label{eqn:fwi}
  J_{\rm FWI}[m;d]= \frac{1}{2}\|F[m]-d\|^2 = \frac{1}{2}\sum_{\bx_s,\bx_r}\int\,dt\,|F[m](\bx_r,t;\bx_s)-d(\bx_r,t;\bx_s)|^2
\end{equation}
The integrals (or discrete approximations) are taken over the recording interval.

The approach to local optimization taken here (and in most work on FWI
and related topics) is based on the gradient of the objective defined
in equation \ref{eqn:fwi}:
\begin{equation}
  \label{eqn:fwigrad}
  g = \nabla  J_{\rm FWI}[m;d] = D_m(F[m])^T(F[m]-d).
\end{equation}
In this formula, $D_m(F[m])$ is the derivative of $F[m]$ with
respect to $m$. 
This gradient is the Euclidean (or $L^2$) gradient, that is, the vector $g$ for
which the Euclidean inner product
\begin{equation}
  \label{eqn:eucip}
  \langle g, \delta m\rangle = g^T\delta m
\end{equation}
with any other vector $\delta m$ gives the
rate of change in the direction of that vector of $J_{\rm FWI}$ at $m$.

The theoretical basis for MSWI mentioned in the Introduction pertains to wave propagation through acoustic materials with model parameter fields (density, sound velocity, etc) which are smoothly varying in space. 
Most of the examples presented below are designed to explore aspects of this theoretical prediction: the search for an optimal model (for either FWI or MSWI) should take place within a set of smooth models. In practice, wave measurements encompass finite ranges of temporal and spatial frequency, and ``smooth'' or ``smoothly varying'' means ``slowly varying on the scale of spatial wavelength''. To ensure that estimated models remain smooth throughout the optimization process, we will penalize oscillation of the
search vector. A {\em weighted inner product} (rather than the
Euclidean inner product $\langle \cdot,\cdot \rangle$) in the space of models, defined by a weight operator $W$, provides a convenient device to apply such a penalty. $W$ should be symmetric and positive definite. As is standard, we define the inner product $\langle \cdot, \cdot\rangle_W$ to be:
\begin{equation}
  \label{eqn:wip}
  \langle g, \delta m\rangle_W = g^TW\delta m
\end{equation}
Comparing the definitions
\ref{eqn:fwigrad}, \ref{eqn:eucip}, and \ref{eqn:wip}, one sees that the
vector $g_W$ for which $\langle g_W, \delta m \rangle_W$ gives the
rate of change of $J_{\rm FWI}$ at $m$ in the direction $\delta m$ is
\begin{equation}
  \label{eqn:fwiwgrad}
  g_W = W^{-1}\nabla  J_{\rm FWI}[m;d] =W^{-1} D_m(F[m])^T(F[m]-d).
\end{equation}
If $W$ is chosen to greatly
amplify components of the vector to which it is applied that oscillate on the wavelength scale and below, 
then those components of $g_W$ must be suppressed relative to the
corresponding components of $g$, hence $g_W$ represents a
non-oscillatory or smooth search direction. That is, the gradient with respect to the weighted inner produce {\em smooths} the gradient \ref{eqn:fwigrad} with respect to the Euclidean inner product. Note that only the inverse operator
$W^{-1}$ appears in the formula \ref{eqn:fwiwgrad}. The weight operator $W$ itself is not needed. Therefore we can conveniently use weight operators $W$ that are difficult to compute, so long as the inverses $W^{-1}$ are easy to compute. 

\subsection{MSWI and Variable Projection}
As mentioned earlier, application of local optimization methods
directly to $J_{\rm FWI}$ tends to produce unsatisfactory model
estimates. MSWI modifies the measure of distance between predicted and
observed data by inserting an adaptive filter field $u$, 
with
one filter per trace. Since both filters and data to which they are applied are defined only on finite time intervals, 
we introduce the {\em truncated filter operator} $K[u]$. 
This operator acts by extending both the filter and the input function be zero outside of their intervals of definition. 
The resulting functions on the real line are convolved then truncated to the domain of the input.
This operator is applied to the predicted data $F[m]$ to
produce the filtered predicted data $K[u]F[m]$.

Note that the filtered predicted data may be viewed as the predicted
data for an {\em extended source}, with the source wavelet at location
$\bx_s$ replaced by $u(\bx_r,\cdot;\bx_s) * w(\cdot)$. That is,
the adaptive filter construction is equivalent to allowing the source
to depend on receiver position as well as source position - one source
wavelet for each source-receiver pair. This construction is an extension of
standard modeling in that the domain (bulk modulus, buoyancy, wavelet)
is larger than in the conventional formulation, and the standard and extended models are the same if the source wavelet is independent of receiver location. This  {\em
  source-receiver extension} \cite[]{HuangSymes2015SEG} is a key ingredient in a number of
other papers on modifications of FWI as mentioned in the Introduction.

It is possible to make the error between filtered predicted data and
observed data as small as one likes by choosing an appropriate filter
field $u$, so this error by itself is useless for estimating the
model. If, on the other hand, $u(\bx_r,t;\bx_s)=\delta(t)$, then the filtered
predicted data is identical to the predicted data.
Therefore, some penalty for divergence of the filter $u$ from
$\delta(t)$ needs to supplement the filtered prediction error. Similar to the works referenced in the Introduction, we use the
mean-square of the filter scaled by a multiple
of time $t$ to measure deviation from $\delta$. To stabilize the computation of the filter we also regularize via a multiple 
of the filter mean-square.
These three terms together give the MSWI objective function:
\begin{equation}
  \label{eqn:filtpen}
  J_{\alpha,\sigma}[m,u;d]=\frac{1}{2}\left(\|K[u]F[m]-d\|^2
  +\alpha^2\|tu\|^2 + \sigma^2\|u\|^2\right).
\end{equation}

The domain space (pairs of bulk modulus fields and adaptive filters)
is very high-dimensional compared with the domain space for FWI (bulk
modulus fields alone, in the present context). It is possible to
optimize $J_{\alpha,\sigma}[m,u;d]$ by alternating updates of $m$ and
$u$  \cite[]{LiAlkhalifah:21}. However the {\em variable projection method} \cite[]{GolubPereyra:73,GolubPereyra:03}
is generally much more efficient than the alternating, or
coordinate search, approach. In this instance, variable projection
consists in minimizing $J_{\alpha,\sigma}[m,u;d]$ over $u$ (a
quadratic optimization) to produce an optimal choice
$u_{\alpha,\sigma}[m:d]$. The reduced objective is
\begin{equation}
  \label{eqn:redfiltpen}
  \tilde{J}_{\alpha,\sigma}[m;d]=\frac{1}{2}\left(\|K[u_{\alpha,\sigma}[m;d]]F[m]-d\|^2
  +\alpha^2\|tu_{\alpha,\sigma}[m;d]\|^2 + \sigma^2\|u_{\alpha,\sigma}[m;d]\|^2\right).
\end{equation}
Note that like the FWI objective, $\tilde{J}_{\alpha,\sigma}[m;d]$
depends only on $m$, with $d$ as a parameter.

\subsection{Relation to Travel Time Tomography}
The first step in relating MSWI to travel time tomography is this simple observation about the small-$\alpha$ behaviour of the reduced MSWI objective $\tilde{J}_{\alpha,\sigma}.$ Namely, that 
\begin{equation}
  \label{eqn:tomo1}
  \lim_{\alpha \rightarrow 0}\frac{1}{\alpha^2} (\tilde{J}_{\alpha,\sigma}[m;d]
  -\tilde{J}_{0,\sigma}[m;d]) = \|t u_{0,\sigma}[m;d]\|^2.
\end{equation}
The right-hand side of equation \ref{eqn:tomo1} is itself a potential objective for a version of MSWI. In fact, it is the variant discussed by \cite{HuangSymes2015SEG,HuangSymes:Geo17}). The related paper by \cite{Symes:24a} establishes the identity \ref{eqn:tomo1}, and explains the relationship between the RHS and the AWI objective function via an asymptotic analysis. 
High-frequency asymptotics are introduced by choosing the source wavelet, $w,$ from a family $\{w_{\lambda}: 0 < \lambda \le 1\}$, defined by
\begin{equation}
\label{eqn:wlam}
    w_{\lambda}(t) = \frac{1}{\sqrt{\lambda}}w_1\left(\frac{t}{\lambda}\right), \, 0<\lambda \le 1.
\end{equation}
The wavelet $w_1$ is assumed to be smooth, of compact support (vanishing for large $t$), and to have zero mean. The scaled wavelets $w_{\lambda}$ inherit these properties, and the dimensionless parameter $\lambda$ 
is proportional to the RMS wavelength of $w_{\lambda}$. Note that the $L^2$ norm of $w_{\lambda}$ is independent of $\lambda$.

Correspondingly, the theoretical basis for MSWI established in \cite{Symes:24a} really pertains to a $\lambda$-indexed family of models and data.
Assume that
\begin{itemize}
    \item[1. ] The regularization parameter $\sigma^2$ appearing in equations \ref{eqn:filtpen}, \ref{eqn:redfiltpen}, and \ref{eqn:tomo1} is proportional to the wavelength proxy parameter $\lambda$ in definition \ref{eqn:wlam} by a factor $r$: 
    \begin{equation}
        \label{eqn:factor}
        \sigma^2 = r\lambda.
    \end{equation}
    The units of $r$ are pressure$^2 \times $ time$^2$, so that definition \ref{eqn:filtpen} is dimensionally consistent. Call the corresponding family of modeling operators $\{F_{\sigma}\}$. These depend implicitly on $\lambda$.
    \item[2. ] Models $m = (\kappa,\beta)$ and $m^* = (\kappa^*,\beta)$ have the single-arrival property. Namely, for any source $\bx_s$ and receiver $\bx_r$, there exists a single ray of geometric optics connecting the source and receiver, with travel times $\tau[m](\bx_s, \bx_r)$ and $\tau[m^*](\bx_s,\bx_r)$ respectively.
    \item[3. ] The data $d = d_{\sigma}$ is model-consistent (or zero noise), that is, $d_{\sigma} = F_{\sigma}[m^*]$.
    \item[4. ] Substitute $F_{\sigma}$ for $F$, $d_{\sigma}$ for $d$ in the definitions \ref{eqn:filtpen} and \ref{eqn:redfiltpen}.
\end{itemize} 

Then for every source-receiver pair $({\bf x}_s,{\bf x}_r)$, there exists $M[m,m^*,w_1]({\bf x}_s,{\bf x}_r) > 0$ so that  
\begin{equation}
    \label{eqn:tomo3}
    \|t u_{0,\sigma}[m;d_{\sigma}](\bx_r,\cdot;\bx_r)\| = \lambda^{-1/2} M[m,m^*,w_1](\bx_s,\bx_r) |\tau[m](\bx_s,\bx_r) - \tau[m^*](\bx_s,\bx_r)| \\ + 
    O(\lambda^{1/2}).
\end{equation}
Thus the traveltime error at $(\bx_s,\bx_r)$ is equal to the $L^2$ norm of the $t$-scaled adaptive filter trace $tu_{0,\sigma}[m;d_{\sigma}](\bx_r,\cdot;\bx_s)$, scaled by $\lambda^{\frac{1}{2}}$ and a $\lambda$-independent factor. The error is proportional to $\lambda$ (an RMS wavelength). From the limit relation \ref{eqn:tomo1} and the inequality \ref{eqn:tomo3},
\[
\lim_{\lambda \rightarrow 0} \lim_{\alpha \rightarrow 0} \frac{\lambda}{\alpha^2}(\tilde{J}_{\alpha,\sigma}[m;d_{\sigma}] -\tilde{J}_{0,\sigma}[m;d_{\sigma}])
\]
\begin{equation}
\label{eqn:tomo4}   
= \sum_{\bx_s,\bx_r} M[m,m^*,w_1]^2(\bx_s,\bx_r) |\tau[m](\bx_s,\bx_r) - \tau[m^*](\bx_s,\bx_r)|^2.
\end{equation}
That is, in the $\lambda \rightarrow 0, \alpha \rightarrow 0$ limit, the reduced MSWI objective scaled by $\lambda/(\alpha^2)$ is proportional to a weighted version of the usual travel time tomography objective 

As shown by \cite{Symes:24a}, $\tilde{J}_{0,\sigma}$ is bounded independently of $\lambda$, so the second term in the left-hand side of equation \ref{eqn:tomo4} can be neglected for small $\lambda$ (small $\sigma$). In the computations carried out below, we presume that for the choices made for $\sigma$, $\alpha$, $\tilde{J}_{\alpha,\sigma}$ can be regarded as a good approximation to the weighted mean square travel time error on the right-hand side of equation \ref{eqn:tomo4}, hence as a suitable proxy for travel time inversion.


This relation is the one mentioned in the introduction. It does not
show that the only stationary points are those of the mean square
travel time error, even approximately, since the proportionality factor $M$ depends on $m$. However if the scale factor $M[m,m^*,w_1]$ is
relatively insensitive to changes in model, as is true if the
source-receiver distance is well away from developing multiple
arrivals, then any stationary point of $\tilde{J}_{\alpha,\sigma}$ is
either close to a stationary point of $\|\tau - \tau^*\|^2$, or far from $m^*$. So we
would expect minimization of $\tilde{J}_{\alpha,\sigma}$ to produce a
model that closely matches the travel times inherent in the data
($\tau[m^*]$ in the notation used here), so long as the models involved are not close to generating multiple arrivals.

AWI adds one more feature, namely scaling {\em each trace} of the filter $u_{0,\sigma}[m;d](\bx_r,\cdot;\bx_s)$ by the reciprocal  of its $L^2$ norm. This trace-by-trace normalization can be interpreted as a choice of weighted
norm on the space of adaptive filters. It leads to a very similar
relation to \ref{eqn:tomo3} for the AWI penalty function, but without
the amplitude factors ($M[m,m^*,w_1]$, scaling by $\lambda$). Thus, the right-hand side of the AWI
analogue of \ref{eqn:tomo4}
is {\em just} the standard travel time tomography objective. See \cite{Symes:24a}
for details.

\cite{Symes:24a} also shows that if energetic multiple arrivals are
present in the data, for any reason, then cross-talk between travel
time branches destroys the relation between MSWI (or AWI) and any
version of the mean-square travel time error. The next sections of
this paper will illustrate the
success of MSWI for single-arrival data and its failure for data with
multiple arrivals.

\subsection{Computational Framework}
The reduced adaptive filter $u_{\alpha,\sigma}[m;d]$ is the solution
of the {\em normal equation}
\begin{equation}
  \label{eqn:normal}
  (S[m]^TS[m] + \alpha^2 t^2 + \sigma^2 I)u = S[m]^Td,
\end{equation}
in which $S[m]u = K[u]F[m]$. This positive definite symmetric linear
system may be solved by various efficient numerical methods. Having
computed $u_{\alpha,\sigma}[m;d]$ (hence the value of 
$\tilde{J}_{\alpha,\sigma}[m;d]$), its (Euclidean) gradient is given by
\begin{equation}
  \label{eqn:gradredfiltpen}
  \nabla \tilde{J}_{\alpha,\sigma}[m;d] =
  D_m(F[m])^TK[u_{\alpha,\sigma}[m;d]]^T(K[u_{\alpha,\sigma}[m;d]]F[m]-d)
\end{equation}
Apart from the appearance of the truncated filter operator
$K[u_{\alpha,\sigma}[m;d]]$, this gradient is almost identical to the FWI
gradient \ref{eqn:fwigrad}. 

Use of a weighted norm in the model space goes exactly as before: with
weight operator $W$, the weighted gradient is
\begin{equation}
  \label{eqn:wgradredfiltpen}
  \nabla_W \tilde{J}_{\alpha,\sigma}[m;d] =
  W^{-1}D_m(F[m])^TK[u_{\alpha,\sigma}[m;d]]^T(K[u_{\alpha,\sigma}[m;d]]F[m]-d).
\end{equation}

The gradient (or weighted gradient), together with the function value, are the
inputs to first-order optimization algorithms such as Steepest Descent and
Limited-Memory Broyden-Fletcher-Goldfarb-Shanno
iterations (\cite{NocedalWright}). Methods more closely related to Newton iteration require
more involved computations (see for instance \cite{Kaufman:75}). 

\section{Numerical Experiments}

\inputdir{project}

This section presents application of MSWI to inversion of 2D synthetic acoustic data. The computational experiments described here illustrate the capabilities and limitations of MSWI. In constructing these experiments we choose length, time, and frequency scales appropriate to crustal seismic exploration. With other scale choices these examples could become 
simplified representations of ultrasound tomography. 

We emphasize that we do not attempt to illustrate the asymptotic relations discussed in the previous section directly. Instead, we make a choice of wavelet $w$ and weight operator $W$ (implicitly, a choice of wavelength and model smoothness) so that the computed data appears to be well-predicted by geometric acoustics, the foundation of the theoretical results reported above. We presume that this choice places our examples far enough into the asymptotic regime described above that the difference between $J_{\alpha,\sigma}$ and the weighted traveltime error $\lambda^{-1} \|M(\tau-\tau^*)\|^2$ is small. Hence minimizing the former produces small enough traveltime error to be a good initial guess for iterative FWI. 

We begin with an example of a low-velocity circular lens anomaly
with a vertical line of sources to the left of the anomaly and a vertical line of receivers to the right. This model produces single arrival data with clearly identifiable first arrivals shifted by well over a wavelength from those produced by a homogeneous bulk modulus model. 
Unsurprisingly, FWI starting from this homogeneous initial guess fails to recover the anomaly. We then use MSWI to generate a bulk modulus estimate. The final adaptive filter corresponding to this bulk modulus estimate is much closer to a band-limited, source-receiver-independent delta than is the adaptive filter for the initial homogeneous model. 
A subsequent FWI iteration beginning with the MSWI-generated model produces a near-perfect data fit as well as an adaptive filter which is a near-perfect band-limited delta.

We next illustrate failure of MSWI experimentally in the case of multi-arrival data. We use an oblate lens model, somewhat more focusing than the circular lens of the first example, and consider two source-receiver geometries. The first acquisition setup is the same as that for the circular lens experiment. In the second oblate lens experiment, the source and receiver lines are moved further apart. FWI starting from a homogeneous initial model fails in both cases. The data produced in the first acquisition case shows signal at multiple arrival times. However, the later arrivals have low enough energy that these multiarrivals are essentially treated as noise by MSWI. Hence MSWI produces a successful starting model for FWI. In the case of the second source-receiver geometry, the larger separation between sources
and receivers leads to a more even division of energy between first and later arrivals, and MSWI fails to estimate a model predicting either set of arrival times, hence fails.

The theoretical justification of MSWI depends on the assumed smoothness of the model generating the observed data to be fit. Indeed, the choice of norm in the space of bulk modulus models forces the MSWI iteration to search through a suite of smooth models. The fourth example explores the effect of non-smoothness in the model generating the observed data. We use a variant of the ``Camembert" model of \cite{GauTarVir:86}, the first published illustration of FWI failure through cycle skipping. 
We make implicit use of the observation of \cite{BarnierBiondi:23a,BarnierBiondi:23b}, that extended inversion methods like MSWI allow independent control of temporal and spatial frequency. We estimate a smooth model via MSWI with model space norm that enforces smoothness, then use the resulting model as a starting point for FWI without enforcing smoothness. The result is an excellent data fit that reveals (some of) the sharp edges in the model.

Finally, all of the experiments described so far are ``inverse crimes''. Namely, the data is noise-free and generated by the same simulation algorithm as used to drive the inversions (both MSWI and FWI). Given the close relation between MSWI and travel time inversion, the stability of the latter \cite[]{Mukhometov:75,TarikereZhou:24} suggests stability of the former. For the fifth example, we add roughly 30\% coherent (model-generated) noise to the data of the first (circular lens) example. We show that FWI initialized on the MSWI-generated model is successful. In fact the algorithm picks out the arrivals corresponding to the noise-free data rather than matching the noise. This final inversion gives a good estimate of the RMS error, which in our usage is a synonym for ``noise".

\subsection{Experimental setup}

The experiments described below share a number of
features which we detail in this subsection. 
Acoustic wave
propagation is simulated via a staggered grid finite difference method
\cite[]{vir86,lev88,Cohen:01} of order 2 in time and 8 in space. The time
step is 
chosen so that the simulation results are
stable, given upper and lower bounds on the parameter fields, difference formulae, and 
spatial sampling. The discrete pressure field is
output over the time range $0 < t < 5$ s with sample interval $0.008$
s at externally specified receiver locations
via piecewise linear interpolation in space and cubic spline
interpolation in time. Isotropic point sources are added
into the acoustic fields at each time step via the adjoints of these
interpolation operations. (See \cite{GeoPros:11} for simulator details.)

The spatial domain for bulk modulus, buoyancy, gradients, etc, is a
rectangular region of 8 km horizontally (``x'') and 4 km vertically
(``z''), sampled on a 20 $\times$ 20 m grid. For all examples except the
oblate lens with larger source-receiver separation and the non-smooth Camembert example, we use 181 receivers spaced 20 m apart, located on a vertical line at
$x$= 5000 m, starting at a
depth $z$ = 200 m. Similarly, twenty isotropic point sources, spaced 150 m apart, occupy a line at $x=3000$ m, starting at a depth of 500
m. Figure \ref{fig:m} shows the sources in green and the receivers in blue.  For the second oblate lens example and the Camembert, source and receiver lines are moved 
further apart, to $x=2000$ m and $x=6000$ m respectively. (Other
parameters remain the same.) Note that this geometry is a 2D
cartoon of that used in cross-well seismic data acquisition. 
The
source wavelets $w$ are all trapezoidal bandpass filtered in the range
$[1.0, 2.5, 7.5, 12.5]$ with a median
frequency of 5.875 Hz corresponding to a median wavelength of
$\approx$ 340 m, and a shortest wavelength of 160 m. The wavelets are centered at
$t=1$ s.

The initial model $m_0$ for MSWI in all of the examples is homogeneous, with $\kappa = $ 4
GPa, $\beta = $ 1 cm$^3$/g. As noted earlier, we regard $\beta$ as
fixed in this series of experiments, so mention of it will be
suppressed. The corresponding data for 20 isotropic point sources at $x=3000$ m and
181 point receivers at $x=5000$ m 
is depicted in
Figure \ref{fig:chwd20}. 
This figure shows
pressure samples ranging between $-0.005$ to $0.005$ GPa. (All similar plots have the same color scale.)

We apply a version of weighted LBFGS optimization \cite[]{NocedalWright} to the
minimization of both $J_{\rm FWI}$ and $\tilde{J}_{\alpha,\sigma}$.
The inverse weight operator ($W^{-1}$
in formula \ref{eqn:fwiwgrad}) is a symmetrized
10-point moving average in both spatial directions, repeated once,
except in the fourth example as discussed below. 
Recall that the weight operator itself is not required. The weighted
gradient is computed via formula \ref{eqn:fwigrad}. The adjoint
derivative $DF[m]^T$ is computed via the {\em adjoint state method}
\cite[]{Chavent:74,GauTarVir:86}, with time reversal of the acoustic fields
implemented via optimal checkpointing
\cite[]{Griewank:92,Griewank:book,Symes:06a-pub}. The LBFGS algorithm uses the weighted
gradient ($g_W$, formula \ref{eqn:fwiwgrad}) to modify the search direction at each iteration,
and the optimal step in this direction is approximated by a simple
backtracking line search algorithm \cite[]{NocedalWright}. The
iteration terminates either when the gradient norm has fallen below
1\% of its initial value or when a maximum number of iterations
is reached.

While the initial model $m_0$ satisfies the upper and lower bound constraints of the staggered grid finite difference algorithm, later model updates may violate these bounds. In fact, after only a few LBFGS iterations, bound violations occur in two of the examples discussed below. The finite difference implementation detects such bound violations and raises an exception, terminating the optimization. In effect, the finite difference parameter bounds add inequality constraints to the optimization of $\tilde{J}_{\alpha,\sigma}$. 
An optimization algorithm formulated for bound constraints, for example LBFGS-B \cite[]{Zhu:97}, could be used to avoid constraint violation and early termination in these examples. 
Since the bounds should be satisfied as strict inequalities at stationary points, we adopt a simpler approach. We reparameterize the wave velocity $c(\bx)$ by a dimensionless field $\gamma(\bx)$, via a logistic function 
\begin{equation}
\label{eqn:logistic}
c(\bx) = a + b \frac{\gamma(\bx)}{\sqrt{1 + \gamma(\bx)^2}}.
\end{equation}
In the relation \ref{eqn:logistic}, the parameters $a$ and $b$ are chosen so that $c(\bx)$ lies strictly between the prescribed velocity bounds input to the staggered grid finite difference algorithm, regardless of the value of $\gamma(\bx)$. The bulk modulus is then computed via the relation $\kappa = \rho c^2$. 
This computation defines a map $\Gamma: \gamma \mapsto \kappa$, a diffeomorphism from the space of grid functions with the discretized $L^2$ norm onto the interior of the feasible set of gridded bulk moduli whose associated velocity fields satisfy the  required bounds. 
In those cases in which LBFGS applied to $\tilde{J}_{\alpha,\sigma}$ results in bound violation, we apply LBFGS to $\tilde{J}_{\alpha,\sigma} \circ \Gamma$ instead. Since $\Gamma$ is a diffeomorphism, stationary points of these two functions are in 1-1 correspondence, however application of LBFGS (or any other local optimization algorithm) to $\tilde{J}_{\alpha,\sigma} \circ \Gamma$ cannot produce a bound violation. Thus this simple change-of-variable device transforms a constrained optimization problem into an unconstrained problem, and is appropriate when the solution is expected to lie in the interior of the feasible set.

The MSWI objective $\tilde{J}_{\alpha,\sigma}$ requires choices of the
parameters $\alpha$ and $\sigma$ and solution of the normal equation
\ref{eqn:normal}. This system is far larger than is convenient to solve
by any variant of Gaussian Elimination, even for our small 2D
examples, so iterative methods are 
necessary. We choose the Conjugate Gradient (CG) method for solution of the normal equations
\cite[]{Dan:67,Steihaug:83,NocedalWright} and stop the
iteration by monitoring the reduction in the normal residual
(the gradient of $J_{\alpha,\sigma}[m,u;d]$ with respect to $u$). The
reduction threshold $\rho$ is thus another necessary input
parameter. For all calculations in this paper, 
we choose $\rho = 0.01$. (Note that unlike $\alpha$ and $\sigma$,
$\rho$ is dimensionless.)

The choice of regularization weights ($\alpha$ and $\sigma$) is a
widely studied topic. We mention some methods for this task in the
Discussion Section. For the set of examples presented here, we take a
simpler approach. 
In the asymptotic theory justifying MSWI, $\sigma>0$ is roughly proportional to $\sqrt{\lambda}$. In this study, we explore neither the choice of $\lambda$ nor the functional dependence of $\sigma$ on $\lambda$. Instead we implicitly select $\lambda$, by choosing a single fixed wavelet, and assign 
$\sigma = 10^{-5}$, a value for which inner CG iterations converge relatively rapidly. {
We choose $\alpha$ by computing the adaptive filter
$u_{\alpha,\sigma}[m_0;d]$ at the homogeneous model $m_0$, for the data $d$ used in the first example below, and for $\alpha=10^k$ for several values of $k$. . The largest such
$\alpha$ that yields an RMS predicted data error
($\|S[m_0]u_{\alpha,\sigma}[m_0;d] - d\|$) less than $0.05 \|d\|$ is
$\alpha = 10^{-4}$. We use
this value of $\alpha$ for all MSWI computations.

\subsection{Recovering a circular lens}

The target bulk modulus field ($\kappa$) for the first example 
is depicted in Figure \ref{fig:m}. It contains an acoustic ``lens''
positioned in the center between $z=$1000--3000 m depth. The background bulk modulus is 4 GPa. At the
center of the lens, the bulk modulus is 2.4 GPa. As
buoyancy ($\beta$) is spatially homogeneous at 1 cm$^3/g$, the background
wave speed is 2000 m/s.

\begin{figure}[htbp]   
\begin{center}
\includegraphics[width=\textwidth]{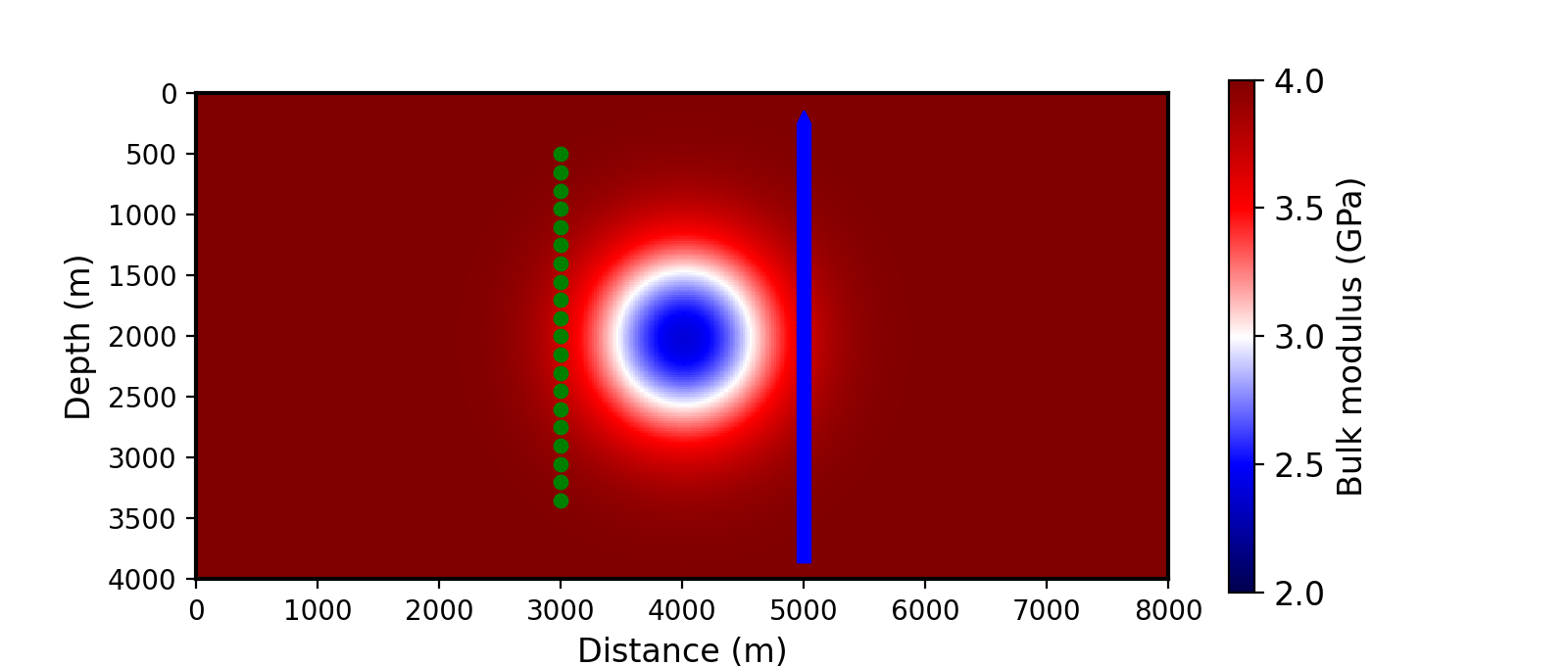}
\caption{Lens model with 20 sources shown in green along a vertical line at $x=3000$ m and 181 receivers shown in blue at $x=5000$ m.}
\label{fig:m}
\end{center}
\end{figure}

\FloatBarrier
\begin{figure}[htbp]   
\begin{center}
\includegraphics[width=\textwidth]{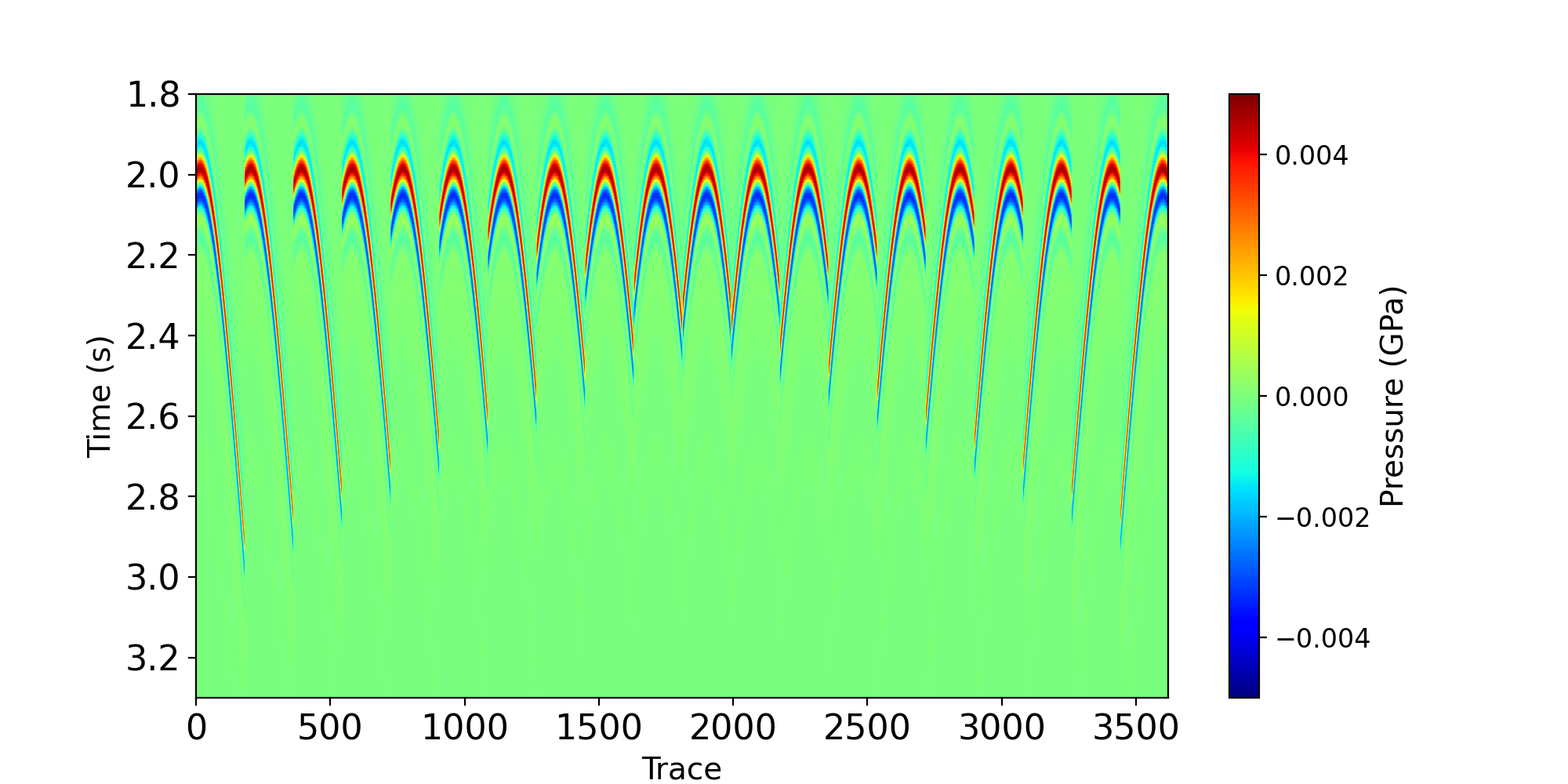}
\caption{Data from the homogeneous initial model $m_0$ acquired using 20 sources located along a vertical line at $x=3000$ m and 181 receivers at $x=5000$ m.}
\label{fig:chwd20}
\end{center}
\end{figure}
\FloatBarrier
The simulated data from this source-receiver configuration for the lens model appears in Figure
\ref{fig:cwd20}.


\begin{figure}[htbp]%
  \centering
  \subfloat[\centering]{
    \includegraphics[width=0.545\textwidth]{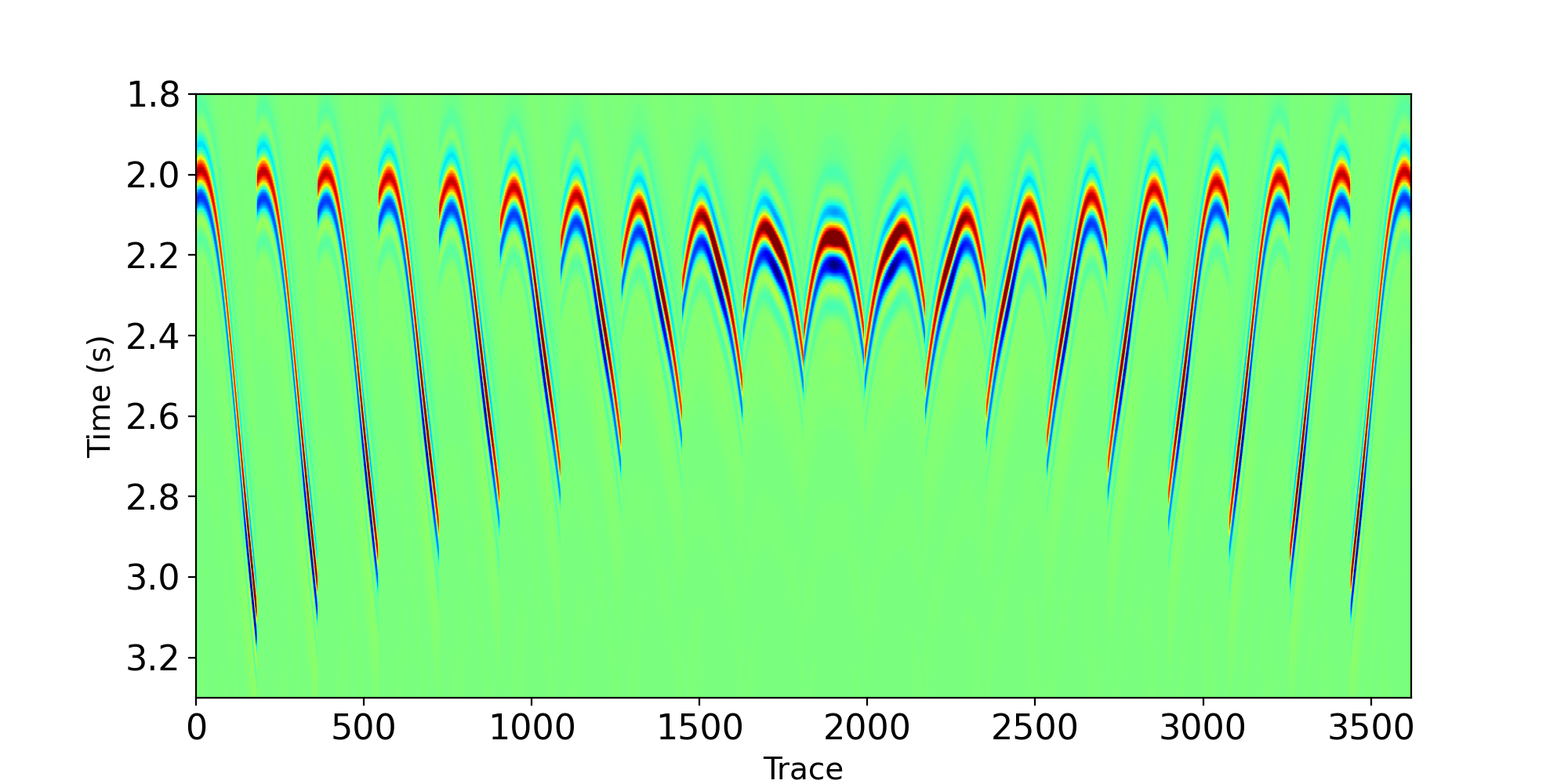}
    \label{fig:cwd20}
  }
     \hspace{-4em}
   \subfloat[\centering]{{\includegraphics[width=0.545\textwidth]{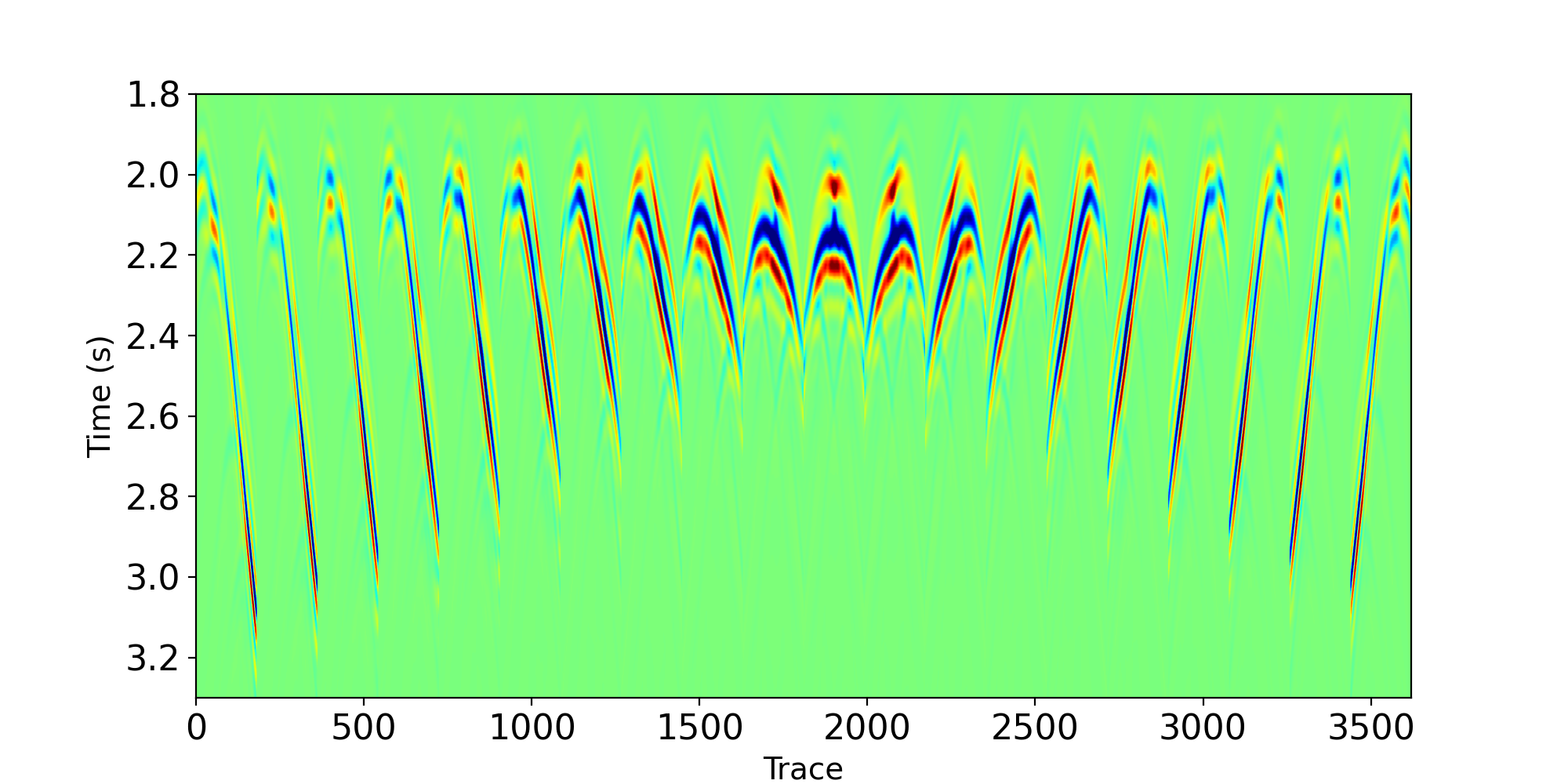} }\label{fig:cwlens20residmestfwi0}}%
   \caption{(a) Data for circular lens model shown in Figure (\ref{fig:m}). (b) FWI residual (difference between
  data predicted from the model shown in Figure \ref{fig:cwlens20mestfwi0} and target data shown in Figure 
  \ref{fig:cwd20}). Note the failure to match arrival times of the later signal, leading to large amplitudes in the
  central part of the display.}%
   \label{cwlens20residmestmswi}%
\end{figure}


The value of $J_{\rm FWI}[m_0,d] \approx 4.6$, and the initial gradient norm at $m_0$ is 2.1 $\times
10^{-3}$. 
After 12 LBFGS steps the objective value has
decreased to 2.4, and the gradient norm to $\approx 3.6 \times 10^{-4}$ (an order of magnitude). The
final FWI inversion result appears in Figure \ref{fig:cwlens20mestfwi0}.


\begin{figure}[htbp]   
\begin{center}
\includegraphics[width=\textwidth]{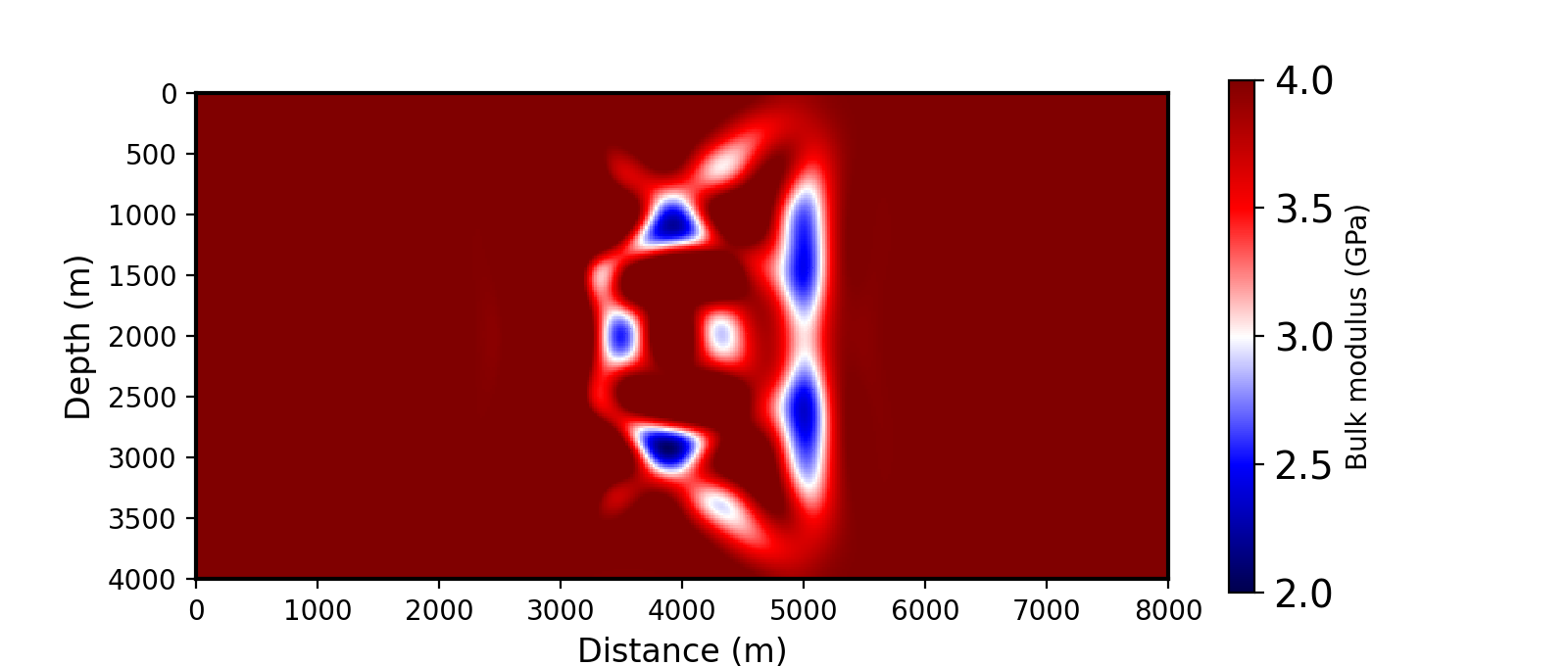}
\caption{Bulk modulus produced by minimizing the FWI objective function given in Equation \ref{eqn:fwi} using 12
  LBFGS iterations.
  The data is shown in Figure \ref{fig:cwd20}.}
\label{fig:cwlens20mestfwi0}
\end{center}
\end{figure}



While the reduction in the gradient norm indicates progress
towards a stationary point, 
the RMS error $\|F[m]-d\|$ is still more than half of its initial value $\|F[m_0]-d\|$. 
Examination of the residual data (Figure \ref{fig:cwlens20residmestfwi0}) shows that the approximate stationary point FWI produces does not predict arrival times correctly. 
In every example shown in this paper, straightforward application of FWI fails to produce a bulk modulus estimate that is even close to accurate. Therefore we will not display any of the FWI bulk modulus estimates or residual data in subsequent examples.
The introduction of the adaptive filter in the formulation of MSWI is intended to overcome this failure to match arrival times. The filter $u$ can shift events to their correct location in time, thus better fitting the data even for drastically incorrect models. The MSWI penalty term in the objective $J_{\alpha,\sigma}[m,u;d]$ measures the extent of the shift. 

The use of a penalty term implies a choice of weight ($\alpha$) as mentioned above. 
This example is the one on which we based the choice $\alpha =
10^{-4}$. We minimized $J_{\alpha,\sigma}[m_0, u; d]$
over $u$ to identify an adaptive filter $u_{\alpha,\sigma}[m_0;d]$. As $\alpha$ increases, the
MSWI fit error
$\|S[m_0]u_{\alpha,\sigma}[m_0;d] - d\|$ at this minimizer increases (and the penalty $\|t u_{\alpha,\sigma}[m_0;d]\|$ decreases).
We choose $\alpha$ just small enough to obtain an RMS error
under 5\% of the data norm $\|d\|$.  The resulting filter $u_{\alpha,\sigma}[m_0;d]$
is shown in Figure \ref{fig:cwlens20uest0}. It exhibits a lot
of energy at non-zero times, as is necessary to move the events in
data computed from the initial model $m_0$ (Figure \ref{fig:chwd20}) to match those in the target data (Figure \ref{fig:cwd20}). 

Throughout the following discussion, we will use the apparent dispersion of energy away from $t=0$ shown by computed adaptive filters to judge the quality of inversions. All plots of adaptive filters in this paper use the same time range $[-0.4,0.4]$ s, amplitude range $[-0.02,0.02]$, and color scale as in Figure \ref{fig:cwlens20uest0}. 



\begin{figure}[htbp]%
   \centering
   \subfloat[\centering ]{{\includegraphics[width=0.5\textwidth]{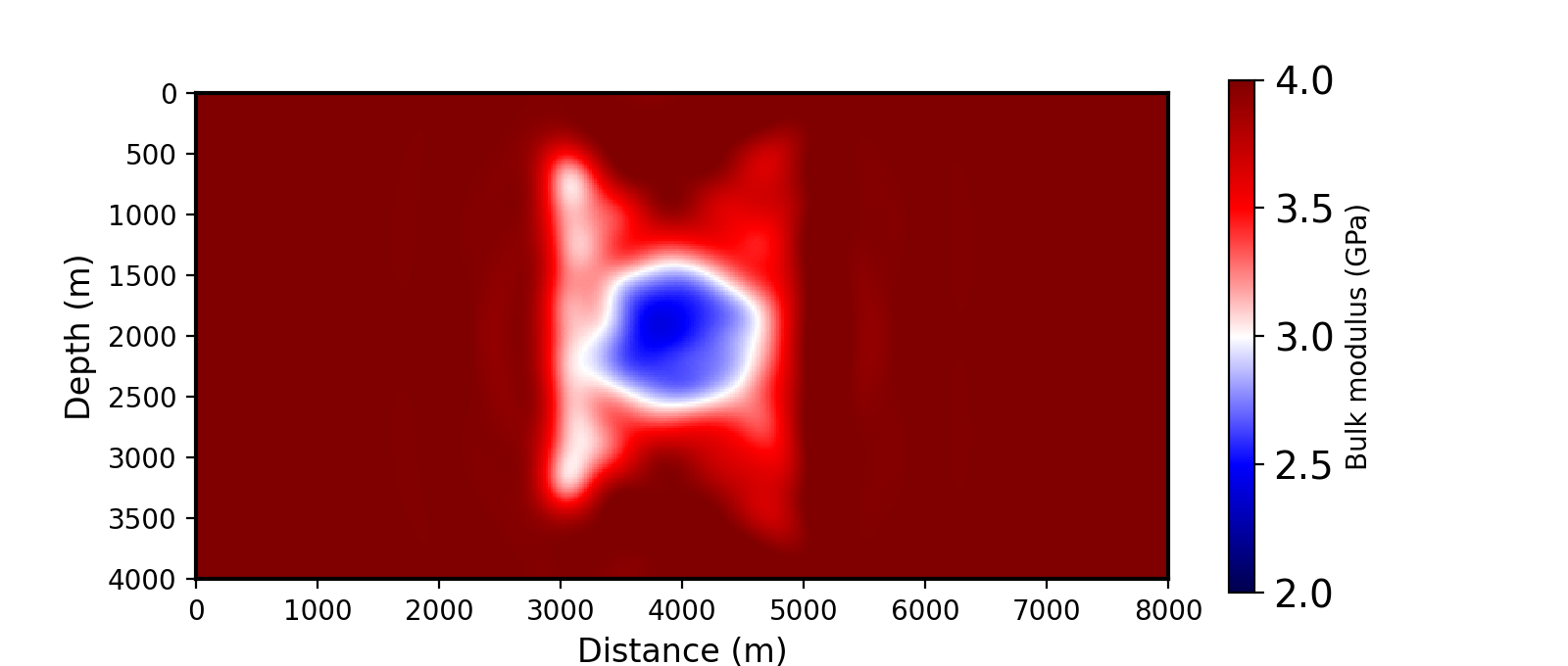}}\label{fig:cwlens20mestmswi}}%
   \hspace{-2em}
   \subfloat[\centering]{{\includegraphics[width=0.5\textwidth]{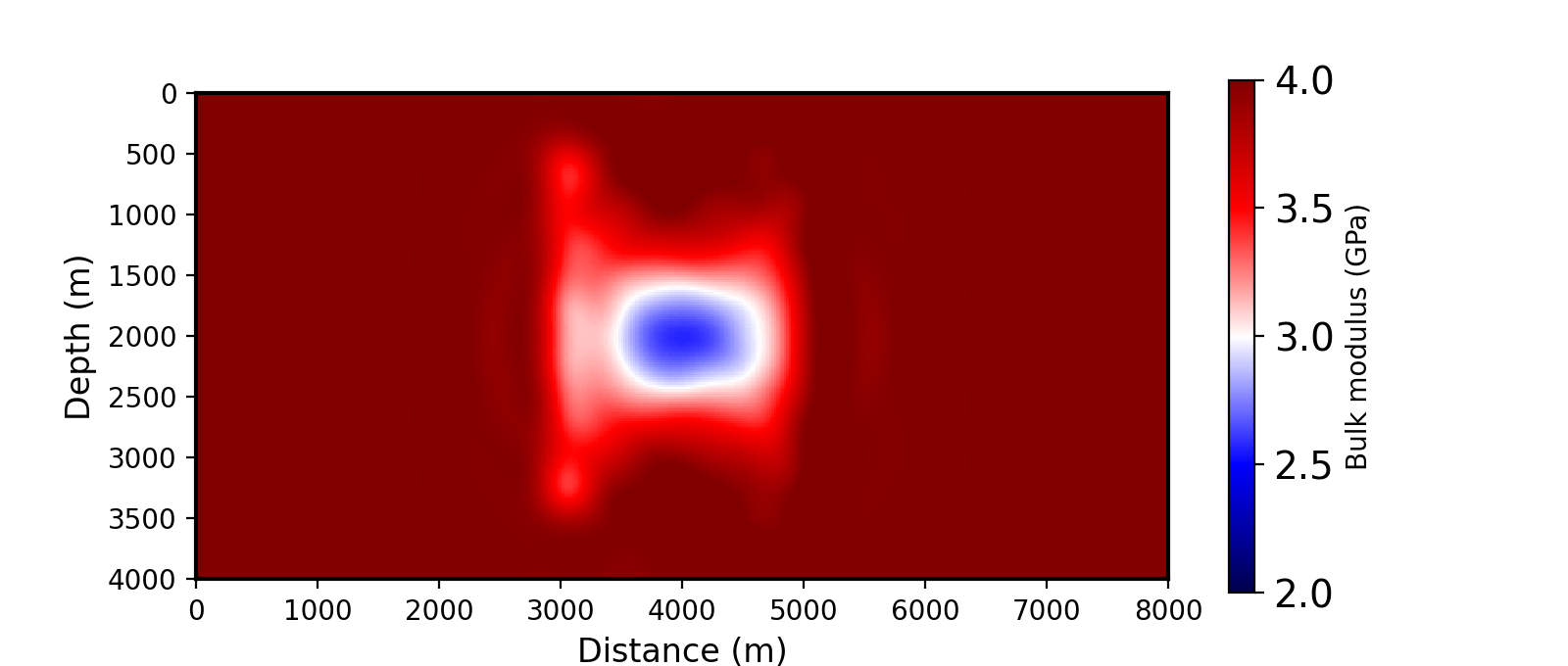} }\label{fig:cwlens20mestmswifwi}}%
   \caption{(a) Bulk modulus produced by 12
  LBFGS steps to minimize the reduced MSWI objective given in Equation \ref{eqn:redfiltpen},
  using the data shown in Figure \ref{fig:cwd20}. (b) Bulk modulus produced by 12
  LBFGS steps to minimize the FWI objective \ref{eqn:fwi},
  using the data shown in \ref{fig:cwd20}, and starting at the MSWI
  result shown in Figure \ref{fig:cwlens20mestmswi}. 
  The gradient norm has been reduced by 2 orders of magnitude.}%
\end{figure}

\begin{figure}[htbp]%
   \centering
   \subfloat[\centering]{{\includegraphics[width=0.5\textwidth]{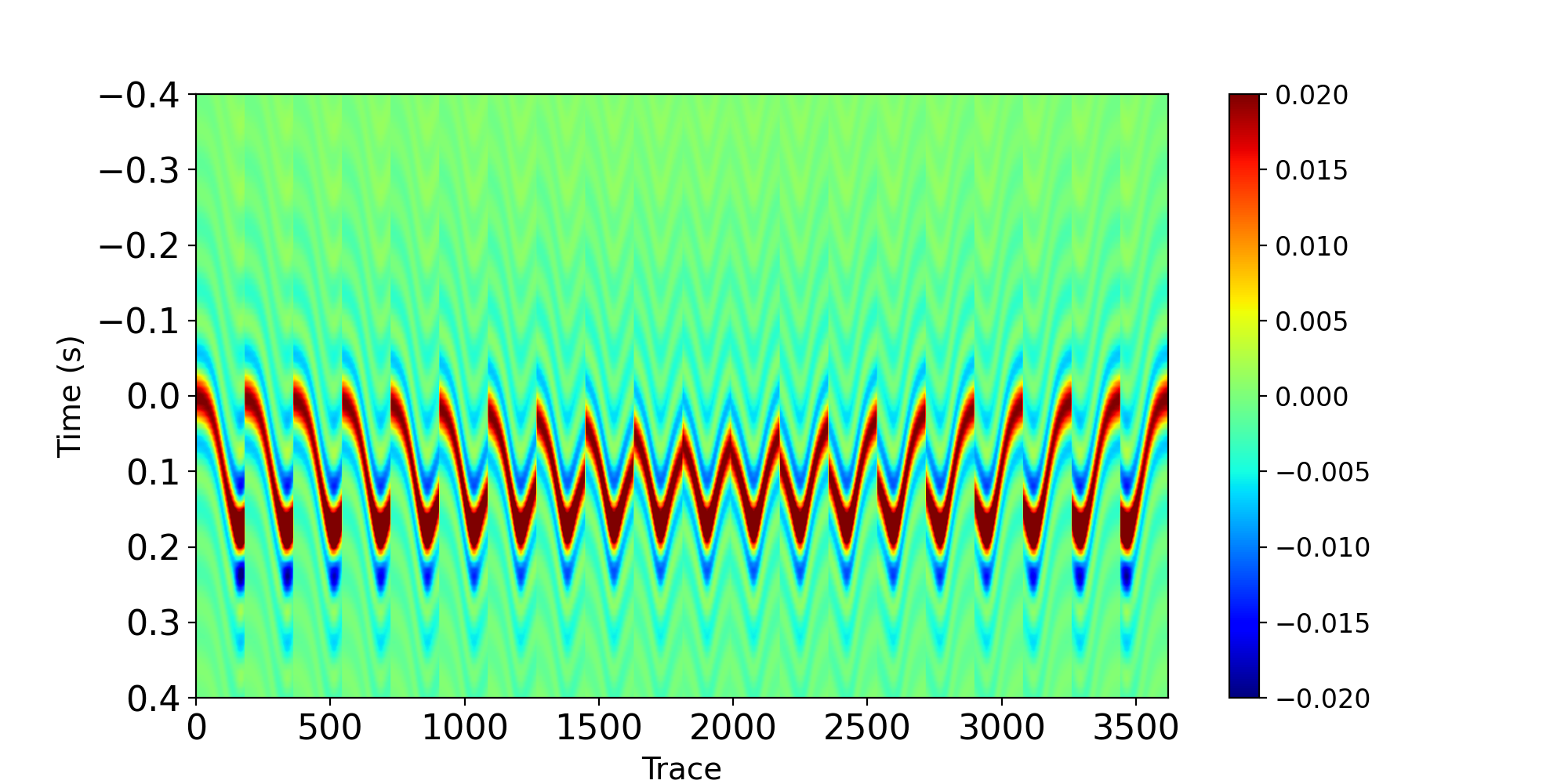}}\label{fig:cwlens20uest0}}%
   \hspace{-2em}
   \subfloat[\centering]{{\includegraphics[width=0.5\textwidth]{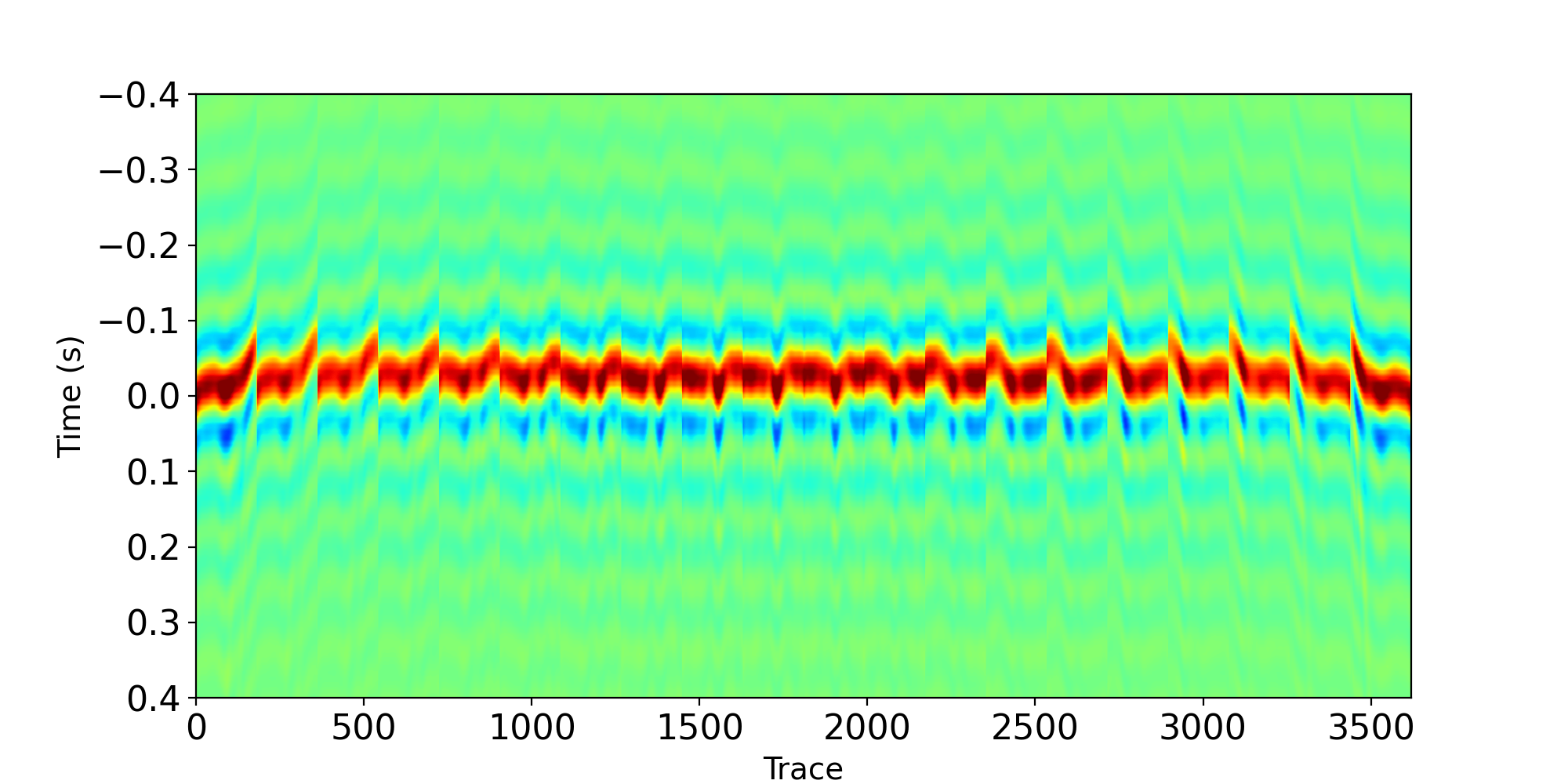}}\label{fig:cwlens20uestmswi}}%
   \caption{(a) Adaptive filter to match
  initial model data (Figure \ref{fig:chwd20}) to circular lens
  data (Figure \ref{fig:cwd20}). Note considerable energy dispersion
  away from $t=0$. (b) Adaptive filter to match
  data computed from the MSWI model (Figure \ref{fig:cwlens20mestmswi}) to circular lens
  data (Figure \ref{fig:cwd20}). Note that there is considerably less energy 
  away from $t=0$ than we see in the filter (Figure \ref{fig:cwlens20uest0}) that matches data from the homogeneous model.}%
\end{figure}


Next we minimize the reduced MSWI objective
$\tilde{J}_{\alpha,\sigma}$ using LBFGS. The initial bulk modulus field is once again
homogeneous at 4 GPa. The initial value of the reduced MSWI objective
is $\approx 1.49 \times 10^{-2}$, and the initial
gradient norm is $\approx 2.2 \times 10^{-5}$. After 12 LBFGS iterations, the objective value has
decreased to $\approx 2.80 \times 10^{-3}$ (a factor of 5.3).
The gradient norm is
$\approx 1.3 \times 10^{-6}$, a reduction of one order of
magnitude. The resulting bulk modulus is shown in Figure
\ref{fig:cwlens20mestmswi}.



The final adaptive filter is shown in Figure
\ref{fig:cwlens20uestmswi}. Note that energy in the filter has
migrated towards $t=0$ compared to the adaptive filter (shown in Figure \ref{fig:cwlens20uest0})
from 
the homogeneous initial model. The majority of
the energy is within an apparent half-wavelength of $t=0$, suggesting that
the MSWI estimate of bulk modulus (Figure \ref{fig:cwlens20mestmswi}) may be an adequate initial estimate for a successful
FWI. 

Starting from the MSWI
inversion result shown in Figure \ref{fig:cwlens20mestmswi}, we minimize the FWI objective function using 12 iterations of LBFGS.
The initial FWI
objective value is $\approx 4.6$, and the initial gradient norm is
$\approx 2.0 \times 10^{-2}$.  After 12 steps of LBFGS,
the prescribed reduction ($10^{-2}$) of the gradient norm is achieved,
and the objective value has decreased to $\approx 2.4 \times
10^{-2}$. The resulting bulk modulus field is depicted in Figure
\ref{fig:cwlens20mestmswifwi}. The data is fit to within roughly 7\% relative RMS error, and the residual (not shown) is almost invisible if plotted on the same scale as the data shown in Figure \ref{fig:cwd20}. 
Thus the combination of MSWI starting from the homogeneous initial model, followed by FWI starting from the MSWI-estimated model, results in a quite precise data fit, far more accurate than that obtained by FWI alone.

\subsection{Oblate lens: single vs. multiple arrivals}

Our next two examples use a more refractive lens model which is no longer circular but instead is extended horizontally, placed somewhat deeper vertically, and has a
minimum bulk modulus of 2 GPa in the center (see Figure \ref{fig:dcwm}.) As is the case with
the circular lens, this oblate lens is smooth 
enough that the predicted data (Figure \ref{fig:dcowd20}) appears to be well-approximated by geometric acoustics.

\begin{figure}[htbp]   
\begin{center}
\includegraphics[width=\textwidth]{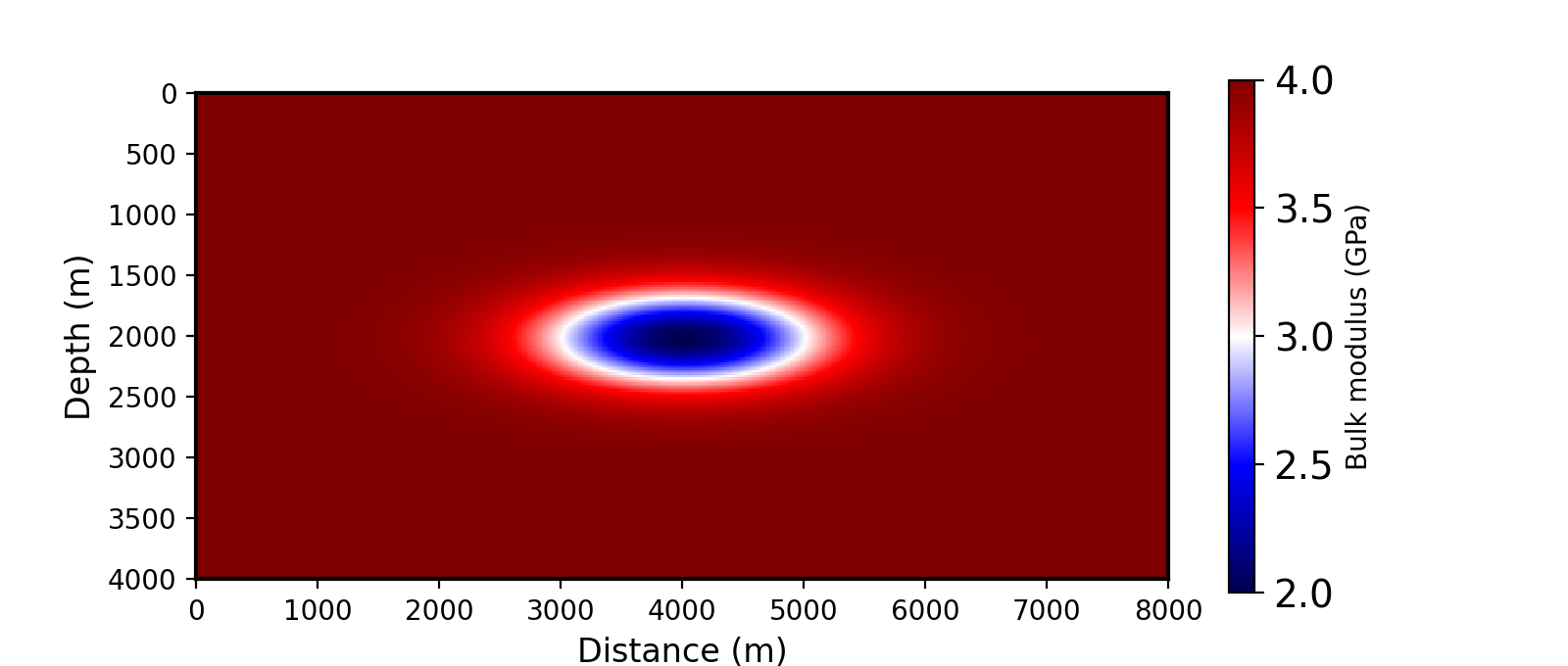}
\caption{Oblate lens model.}
\label{fig:dcwm}
\end{center}
\end{figure}

The first example using this model has exactly the same acquisition
geometry as was used in the circular lens example, namely, 20 sources in a 
vertical line at $x$ =
3000 m, and 181 receivers at $x$ = 5000 m with the same spacing as described in the
previous section. The data is shown in Figure \ref{fig:dcowd20}. 
The stronger refraction of this model results in a visible triplication at
the center of the source-receiver line, and in fact in all parts of
the data, although the smearing effect of finite bandwidth obscures this
detail in several places. Hence most of the data can be (just!) regarded as exhibiting a single
arrival on an identifiable wavefront.

\begin{figure}[htbp]   
\begin{center}
\includegraphics[width=\textwidth]{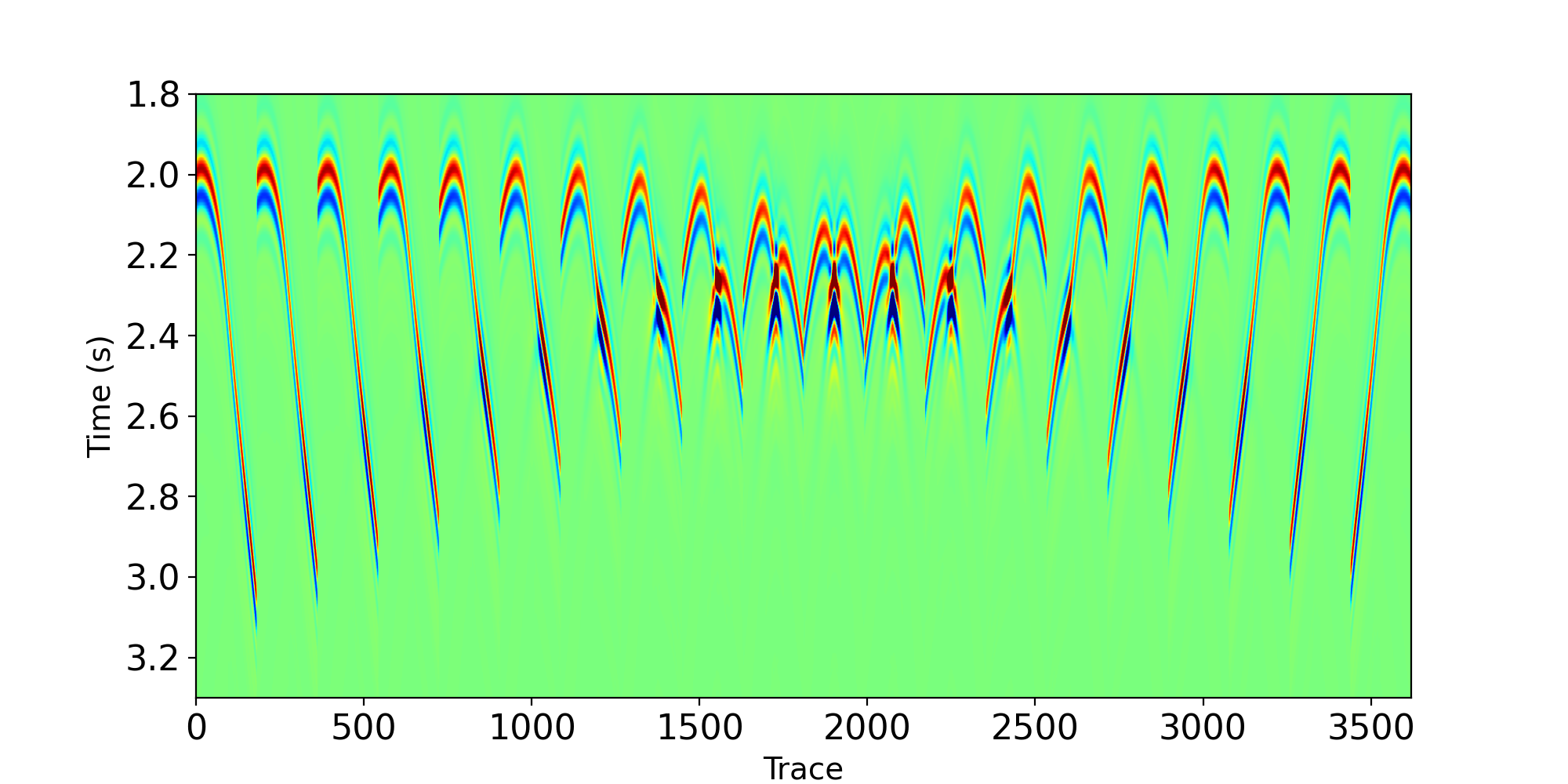}
\caption{Data generated from the model shown in Figure
  \ref{fig:dcwm} and the same source-receiver geometry as for the circular lens
  data.}
\label{fig:dcowd20}
\end{center}
\end{figure}

As already mentioned, FWI initiated from the homogeneous
model fails for all examples shown in this paper. 
Application of LBFGS to the MSWI objective (with the same
parameters as in the previous example) results in the model depicted
in Figure \ref{fig:mnqcaustic20mestmswi} and reduces the objective value and gradient to
roughly 18\% and 7\% of their initial values respectively. 


\begin{figure}[htbp]%
   \centering
   \subfloat[\centering ]{{\includegraphics[width=0.5\textwidth]{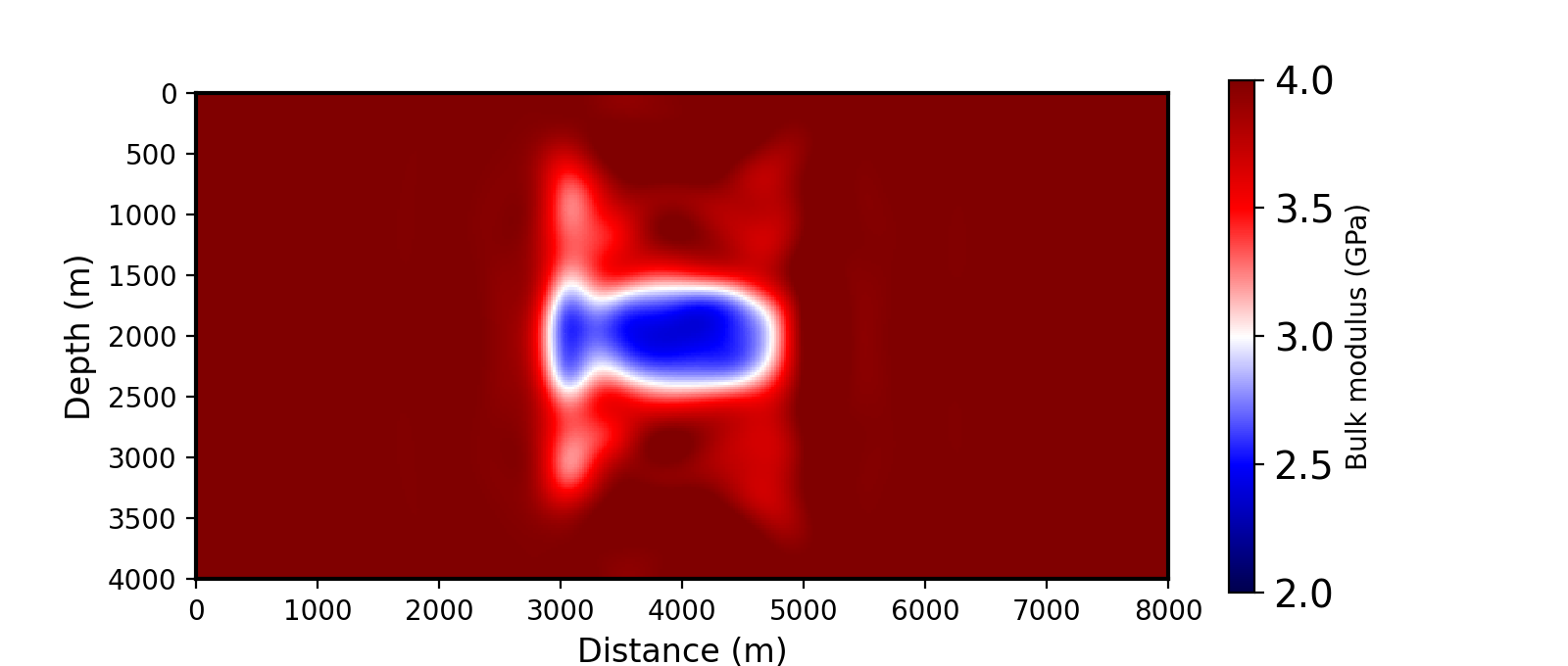}}\label{fig:mnqcaustic20mestmswi}}%
   \hspace{-2em}
   \subfloat[\centering]{{\includegraphics[width=0.5\textwidth]{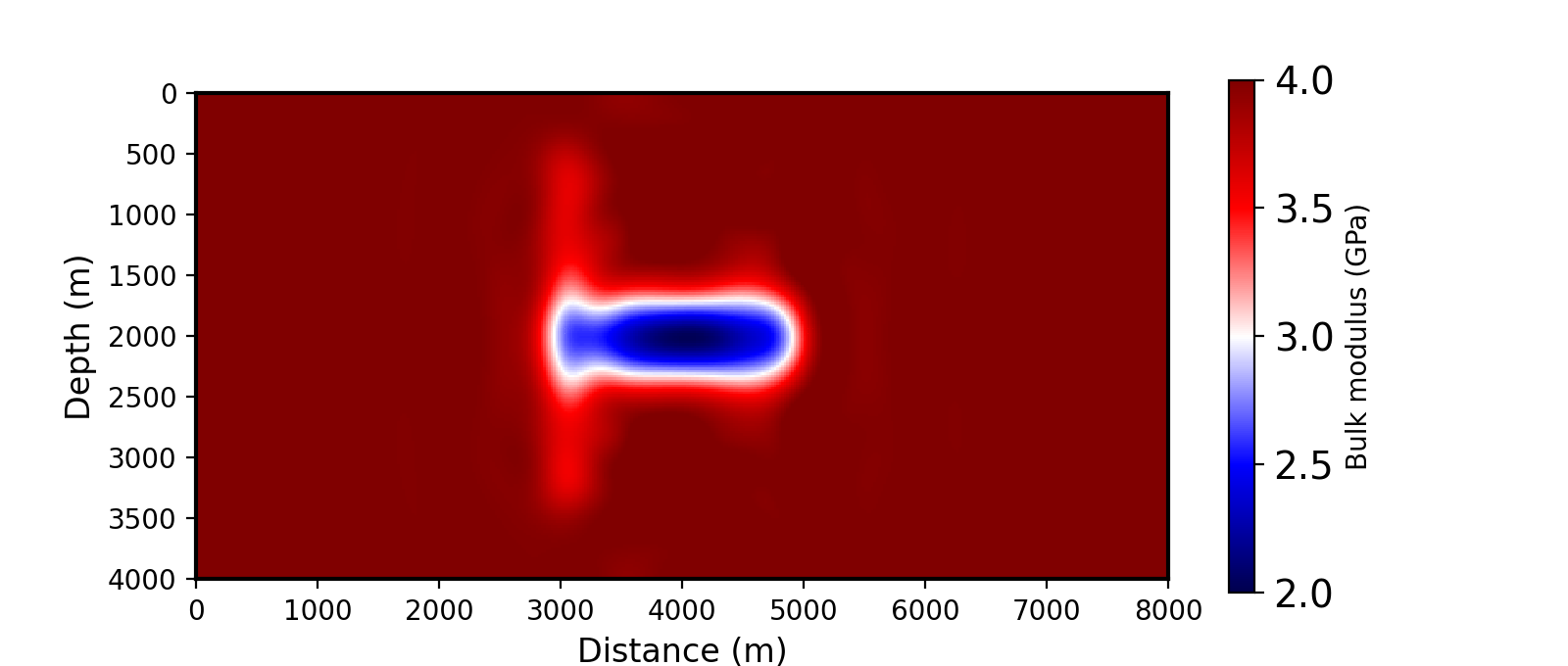} }\label{fig:mnqcaustic20mestmswifwi}}%
   \caption{(a) Result of 12 LBFGS 
  iterations applied to the MSWI objective function and the data shown in Figure \ref{fig:dcowd20} starting from 
  the homogeneous model. (b) Bulk modulus produced by 12
  LBFGS steps applied to the FWI objective 
  using the data shown in \ref{fig:dcowd20} and starting at the approximate MSWI minimizer shown in Figure \ref{fig:mnqcaustic20mestmswi}.}%
\end{figure}


The change in adaptive filters between the homogenous initial model
and the MSWI inversion is shown in Figures
\ref{fig:mnqcaustic20uest0} and \ref{fig:mnqcaustic20uestmswi}. The
reduction in energy dispersion is evident.
The adaptive filter appears to have the bulk of its energy within a half-wavelength of $t=0$, making the MSWI inversion result a plausible initial
estimate for FWI. We applied twelve LBFGS iterations to the FWI objective
function using the data shown in Figure \ref{fig:dcowd20}. The inversion starts from the final MSWI model shown in 
Figure \ref{fig:mnqcaustic20mestmswi}. Figure \ref{fig:mnqcaustic20mestmswifwi}
shows the final model recovered by FWI. The
objective function decreased to 1\% of its initial value.
Thus in this example, as in the first example, iterative LBFGS minimization of the MSWI objective with a modest number of iterations produces a satisfactory initial model for FWI inversion. 


\begin{figure}[htbp]%
   \centering
   \subfloat[\centering ]{{\includegraphics[width=0.545\textwidth]{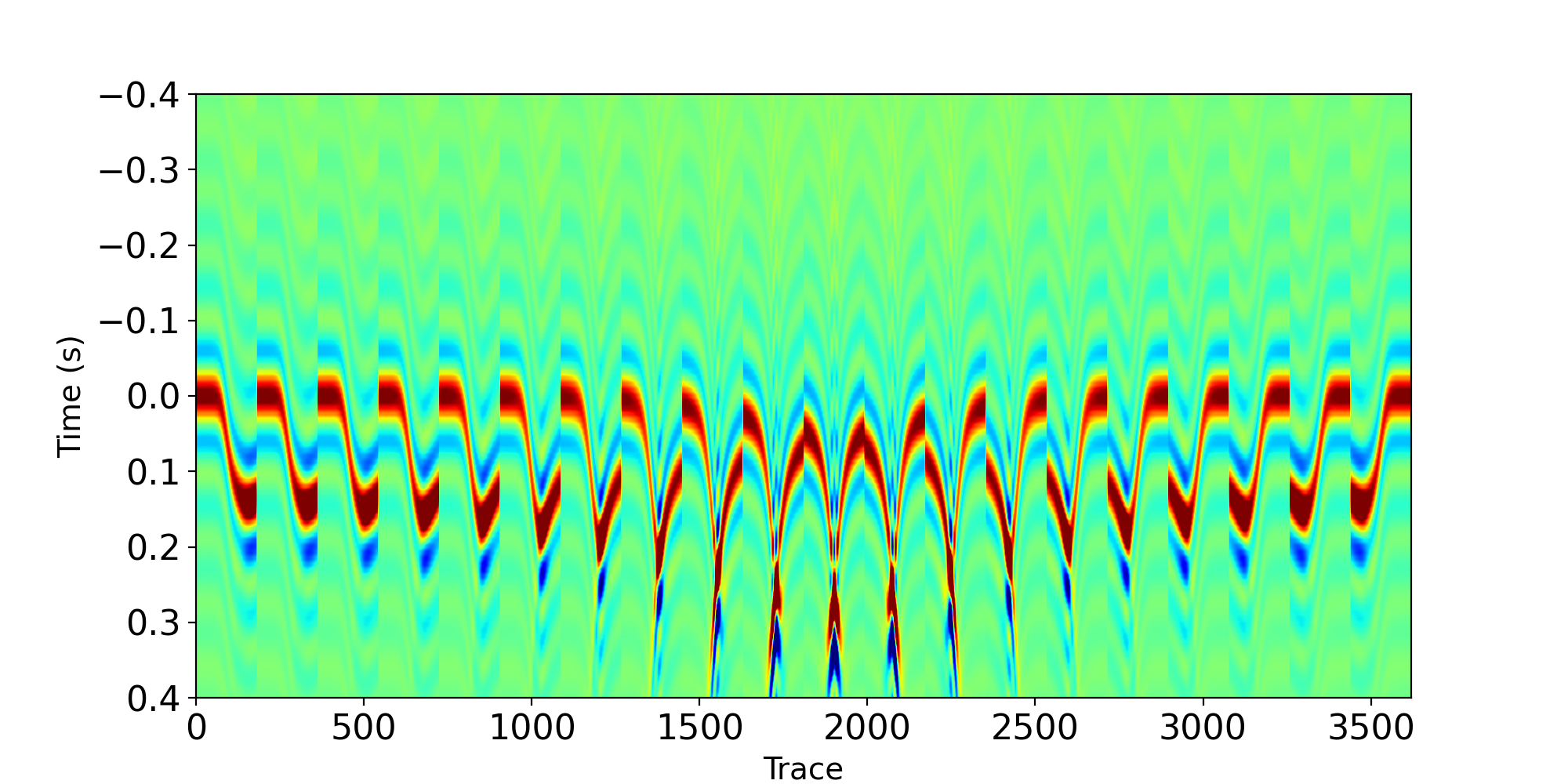}}\label{fig:mnqcaustic20uest0}}%
     \hspace{-4em}
   \subfloat[\centering]{{\includegraphics[width=0.545\textwidth]{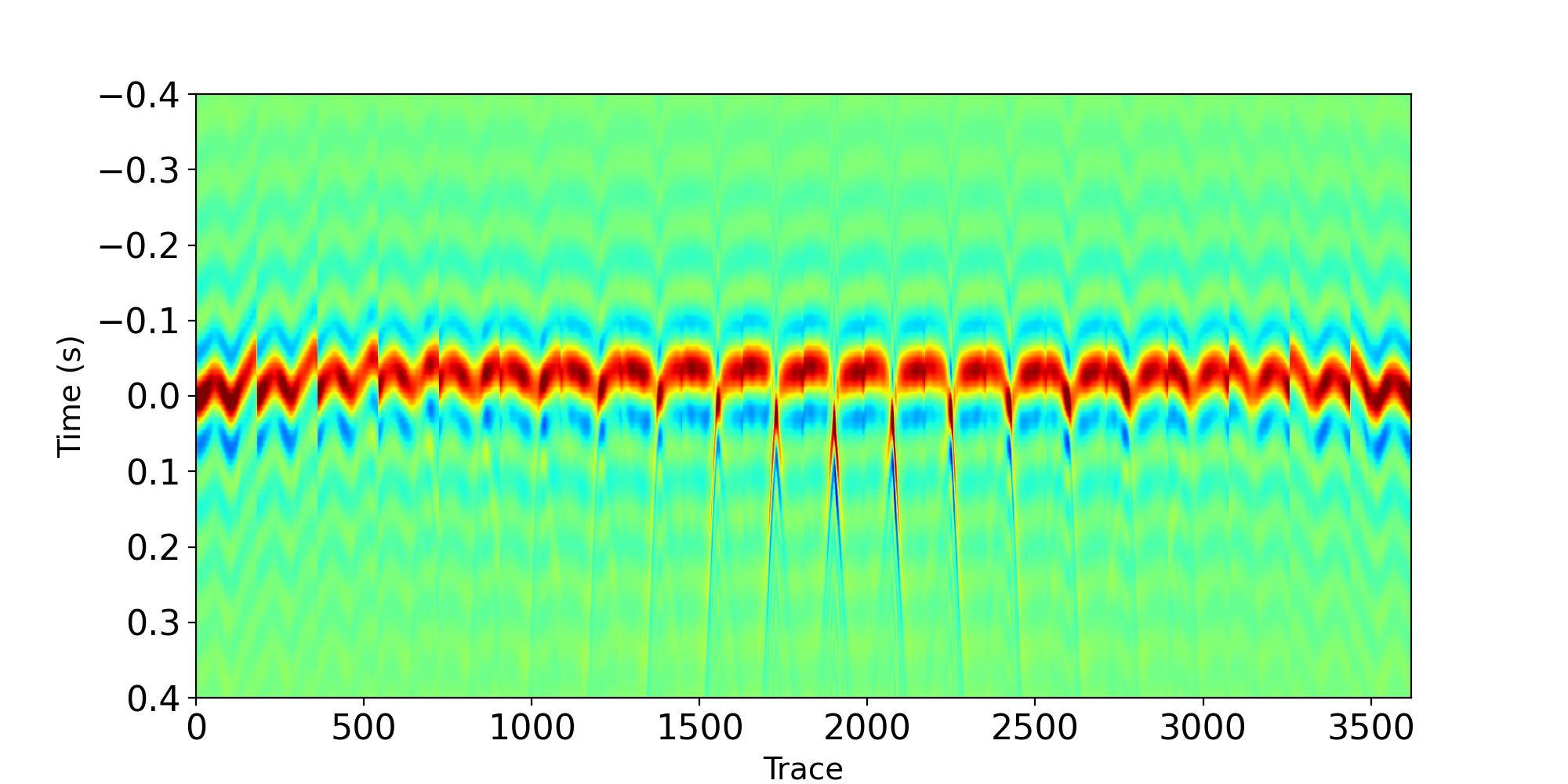} }\label{fig:mnqcaustic20uestmswi}}%
   \caption{(a)
  Adaptive filter estimate starting the inversion at a homogeneous initial model and fitting the data in
  Figure \ref{fig:dcowd20}. (b) Adaptive filter estimate from MSWI inversion result shown in Figure \ref{fig:mnqcaustic20mestmswi}.}%
\end{figure}

The second example based on the target model in Figure \ref{fig:dcwm}
differs from the first {\em only} in that the sources and receivers
have been moved 2 km further apart. The sources now lie on the
line $x=2000$ m, and the receivers are at $x=6000$ m. 
All other aspects of data generation are
the same. The resulting data is shown in Figure \ref{fig:dcwd20}.
Notice that this increase in offset between sources and receivers has allowed distinct energetic later
arrivals to appear in the data.
\begin{figure}[htbp]   
\begin{center}
\includegraphics[width=\textwidth]{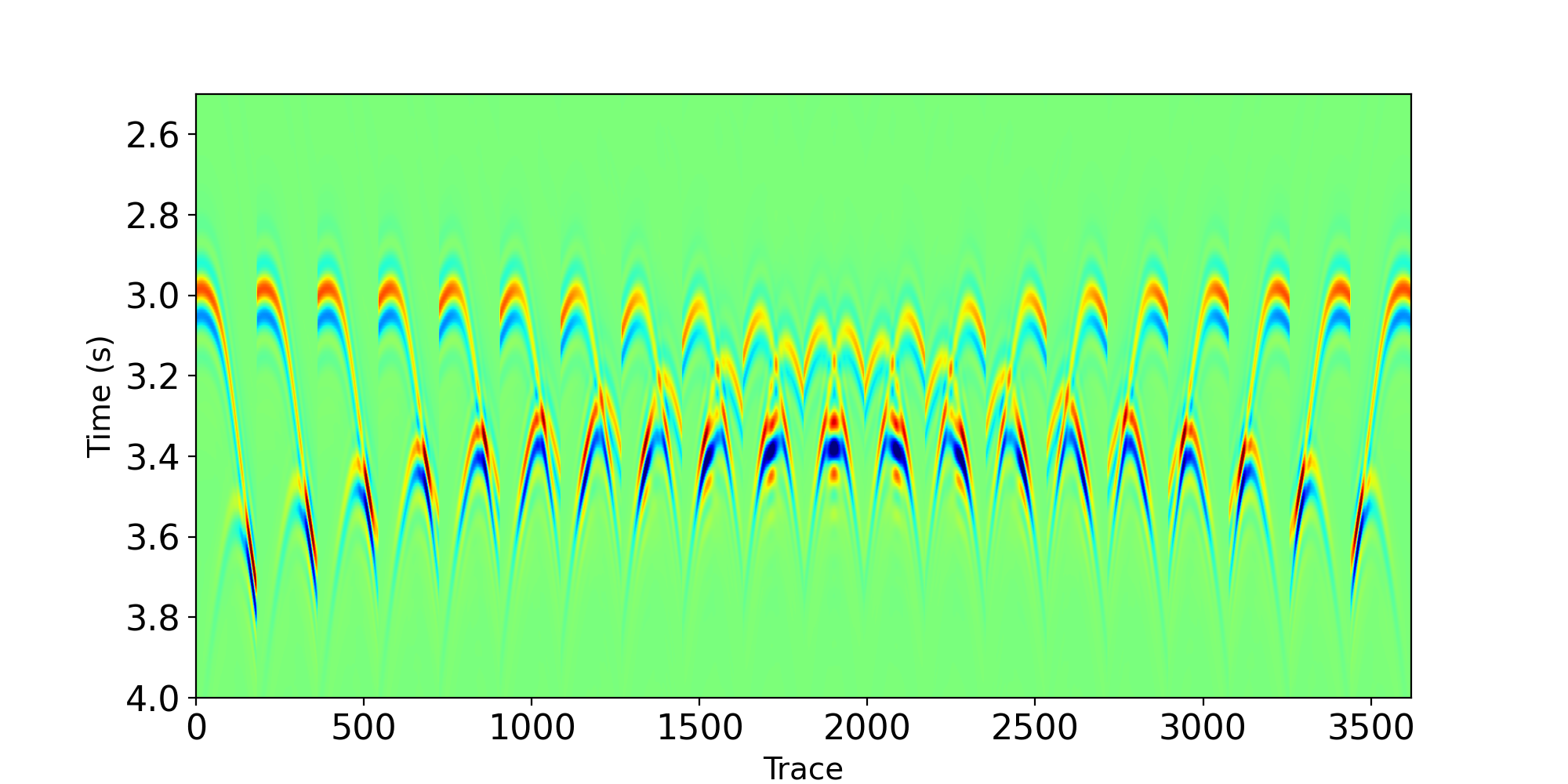}
\caption{Data generated using the same inputs
  as in Figure \ref{fig:dcowd20}, but with sources moved to $x=2000$ m
  and receivers to $x=6000$ m.}
\label{fig:dcwd20}
\end{center}
\end{figure}
As before, we compute the adaptive filter required to produce the data
of Figure \ref{fig:dcwd20} from the corresponding homogeneous medium data (Figure
\ref{fig:chwd20}). This filter, displayed in \ref{fig:mcaustic20covuest0},
shows that for most traces, two or more major energy peaks appear in the
adaptive filter, forced by the need to match multiple energy peaks in
the target data with the single energy peak in the homogeneous medium
data.

\begin{figure}[htbp]   
\begin{center}
\includegraphics[width=\textwidth]{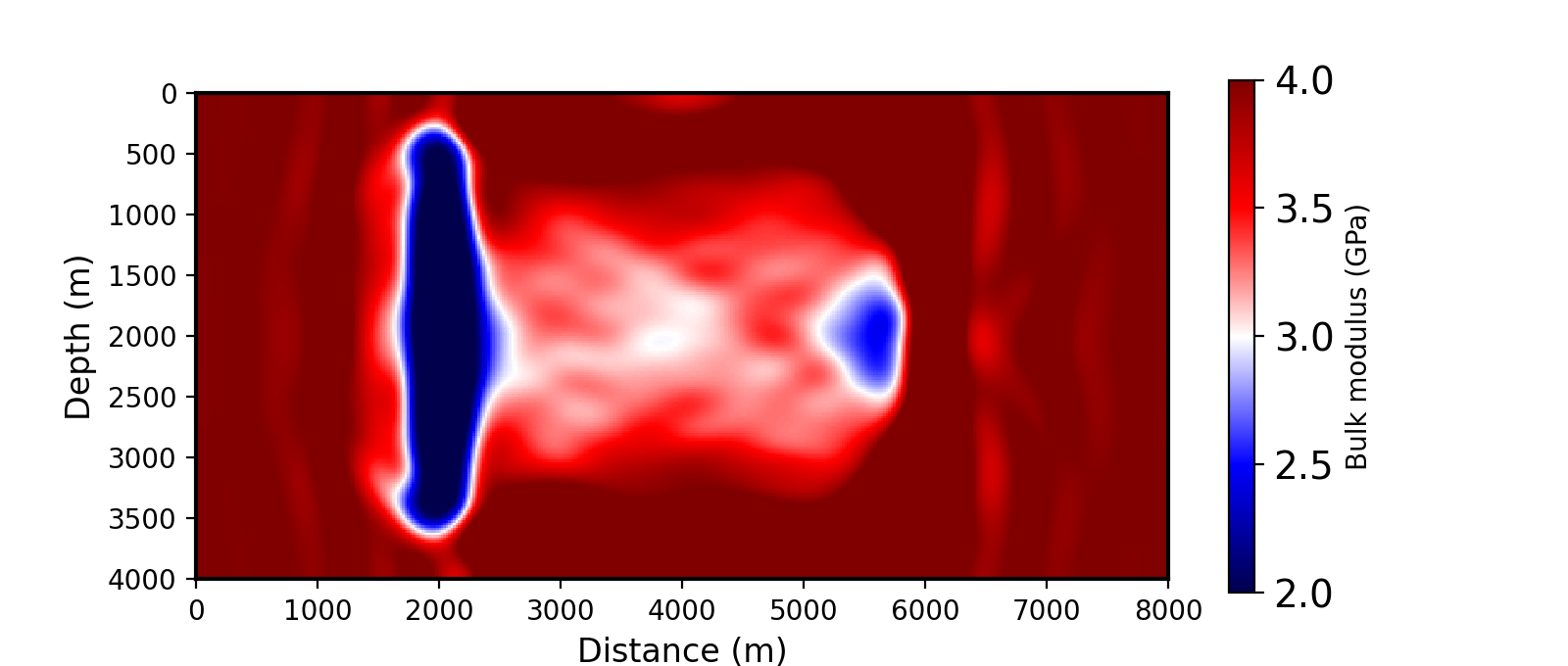}
\caption{Result of 37 LBFGS 
  iterations applied to the  MSWI objective function for the data in Figure \ref{fig:dcwd20} starting from
  the homogeneous initial model.}
\label{fig:mcaustic20covmestmswicont}
\end{center}
\end{figure}

\begin{figure}[htbp]%
   \centering
   \subfloat[\centering ]{{\includegraphics[width=0.545\textwidth]{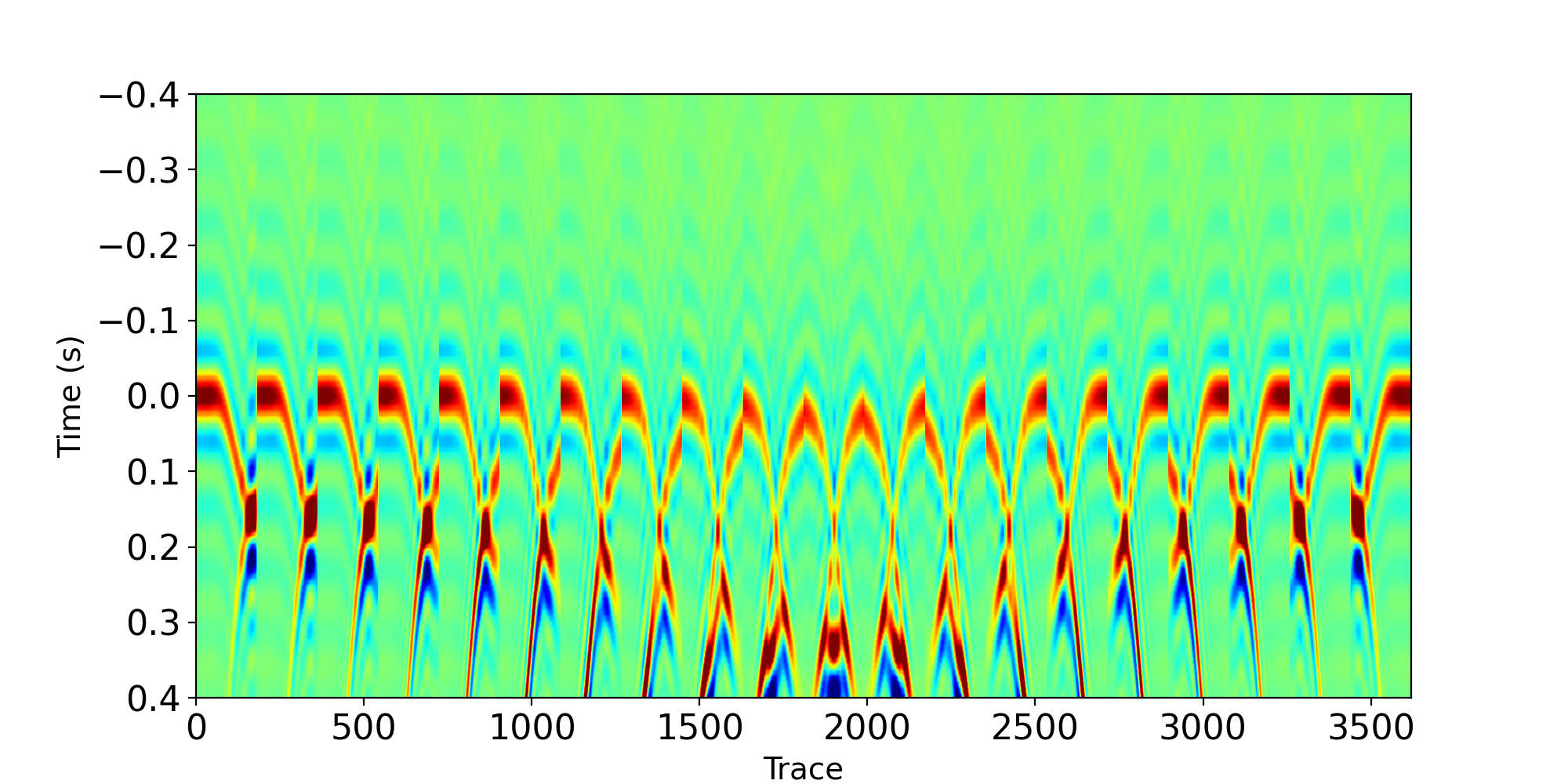}}\label{fig:mcaustic20covuest0}}%
   \hspace{-4em}
   \subfloat[\centering]{{\includegraphics[width=0.545\textwidth]{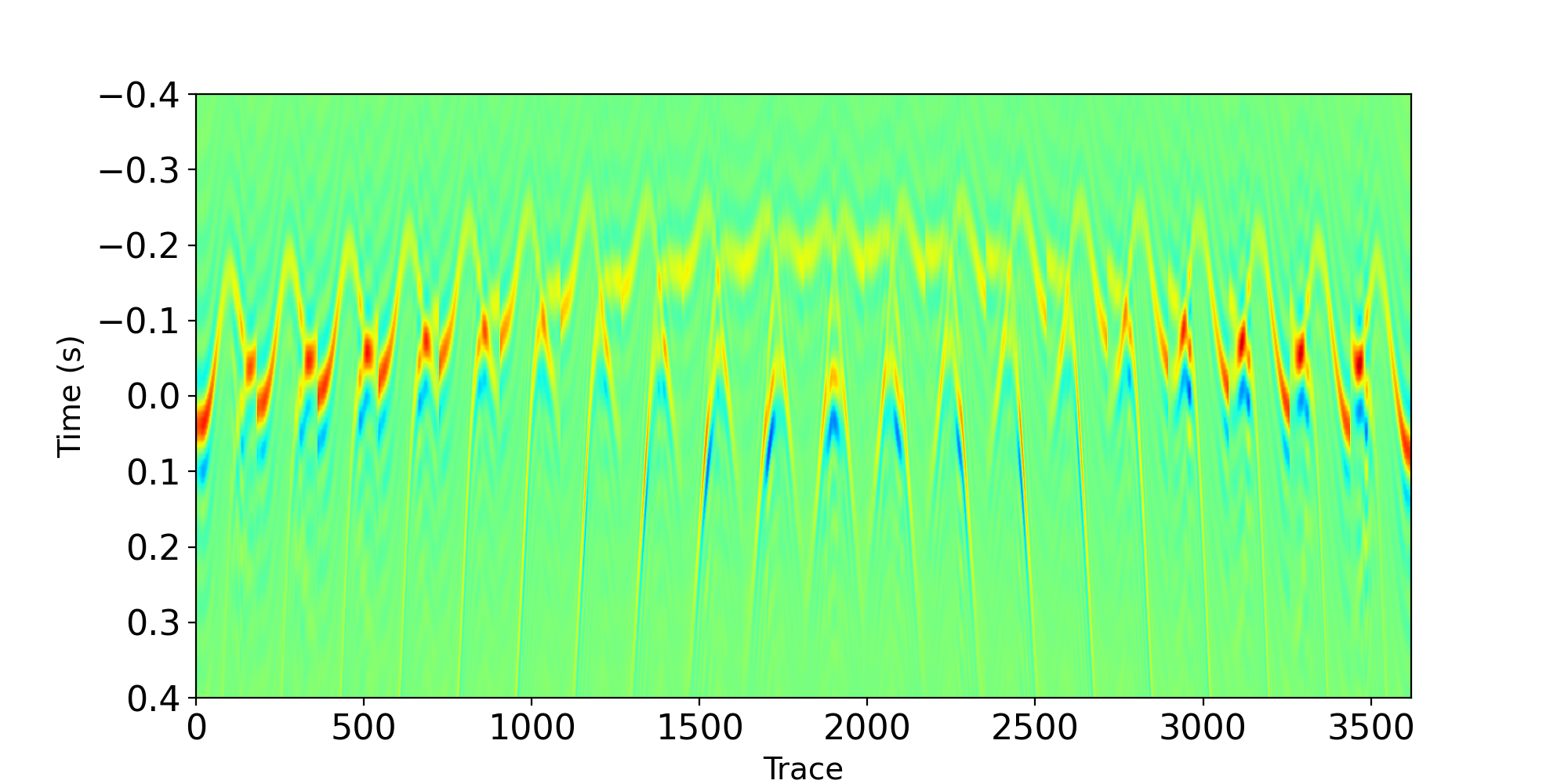} }\label{fig:mcaustic20covuestmswicont}}%
   \caption{(a) Adaptive filter required to produce the data
of Figure \ref{fig:dcwd20} from corresponding homogeneous medium data (Figure
\ref{fig:chwd20}). (b) Adaptive filter estimate for MSWI
  inversion shown in Figure \ref{fig:mcaustic20covmestmswicont}.}%
\end{figure}

This optimization appears to be more difficult. Because the reduction in MSWI gradient and objective function values achieved in 12 LBFGS iterations is not as great as was the case in
the two previous examples, we perform 25 additional LBFGS iterations. 
The final objective value obtained after 37 LBFGS steps is roughly 7\% of the initial objective value, and the gradient norm is approximately 4\% of its initial value. While these norm decreases suggest a reasonably successful optimization, the approximate minimizer (Figure
\ref{fig:mcaustic20covmestmswicont}) does not at all resemble the target bulk modulus (Figure \ref{fig:dcwm}).
The adaptive filter derived as a by-product of the MSWI minimization (Figure \ref{fig:mcaustic20covuestmswicont}) is no more
focused at $t=0$ than is the filter for homogeneous medium data (Figure \ref{fig:mcaustic20covuest0}).


The model shown in Figure \ref{fig:mcaustic20covmestmswicont} does not appear to be a good initial estimate for FWI, and indeed it is not. Applying 25 LBFGS iterations to the FWI objective function starting from the MSWI model produces a reduction of the gradient norm to 3.5\% of its initial value, thus an approximate stationary point. However the value of the RMS fit error has only decreased by a couple of percent. Since the data fit remains unsatisfactory, it is not surprising that the resulting bulk modulus estimate (not shown) 
is no better than the MSWI estimate (Figure \ref{fig:mcaustic20covmestmswicont}).

This pair of examples based on the oblate lens model 
illustrate the behaviour of MSWI described in the {\em Theory Section}.
The data in the first example (Figure \ref{fig:dcowd20}) shows clear signs of caustic formation, but the amount of energy in later arrivals is small. Thus MSWI treats the data as single arrivals plus noise and delivers a model that largely matches the dominant first arrival times. The second example in the pair (see the data in Figure \ref{fig:dcwd20}) distributes energy much more evenly between first and later arrivals, and MSWI fails to produce a useful model estimate. 


\subsection{Recovering a Camembert: effect of non-smoothness}

The justification for MSWI relies on geometric asymptotics. 
The previous examples involve the recovery of a smooth bulk modulus model. To explore the implications of model non-smoothness - a common aspect of prototypical mechanical parameter distributions in the earth, the human body, and other potential wave imaging arenas - we apply the methodology developed above to a variant of a model introduced in one of the the first important papers on FWI, \cite{GauTarVir:86}. This model, shown in Figure \ref{fig:cambulk}, 
features a circular inclusion inside of which the bulk modulus is 20\% higher than in the surrounding. (From its appearance in perspective plot, this model has come to be called the Camembert.)
The geometry of sources and receivers is the same as in last (third) example. The corresponding data (Figure \ref{fig:camcwd20}) differs from that shown in Figure \ref{fig:dcowd20} in that only a single arrival is apparent. Indeed, the Camembert is a defocusing lens, rather than a focusing lens as in the previous examples, and does not generate multiple arrival times.

\begin{figure}[htbp]   
\begin{center}
\includegraphics[width=\textwidth]{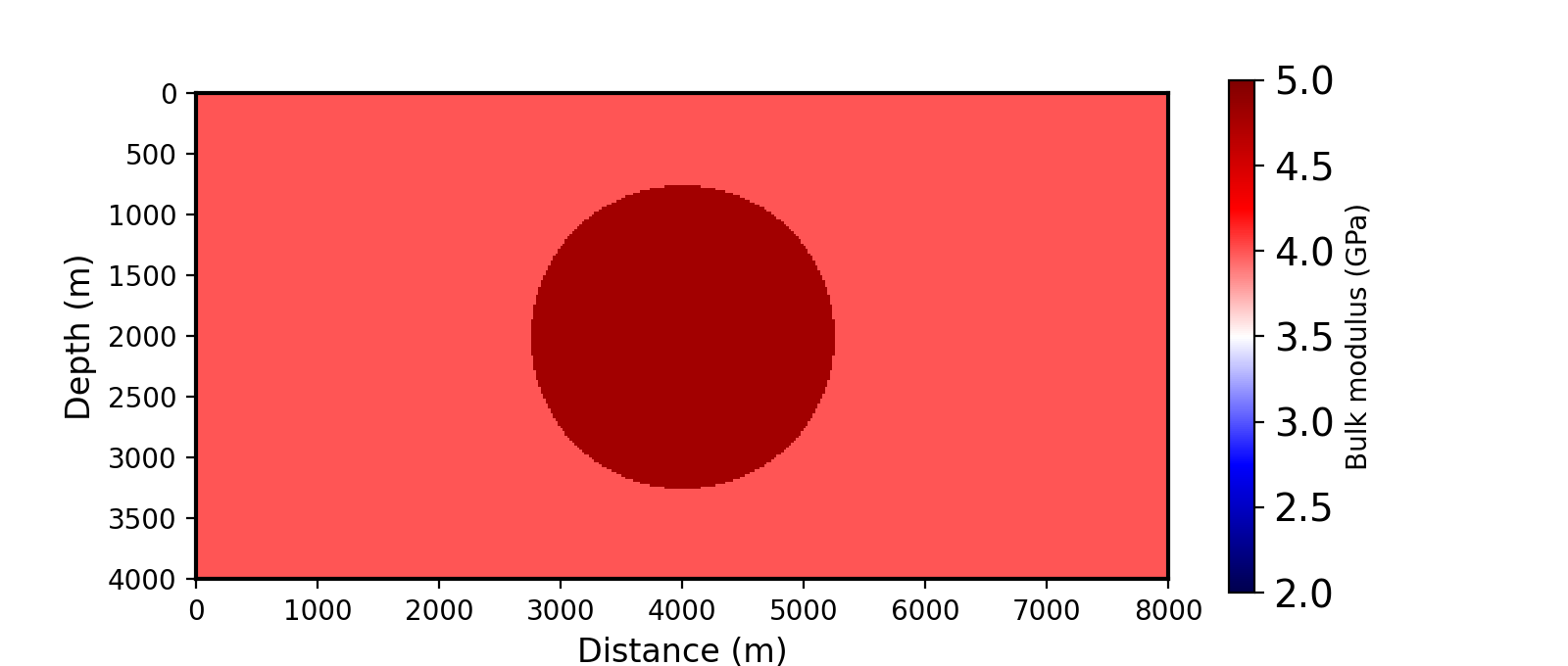}
\caption{A model resembling the Camembert
of \cite{GauTarVir:86}. Diameter of circular inclusion = 2500 m,
approximately 8 wavelengths at median frequency. Interior value (4.8 GPa) is 20\%
above exterior (4.0 GPa).}
\label{fig:cambulk}
\end{center}
\end{figure}

\begin{figure}[htbp]   
\begin{center}
\includegraphics[width=\textwidth]{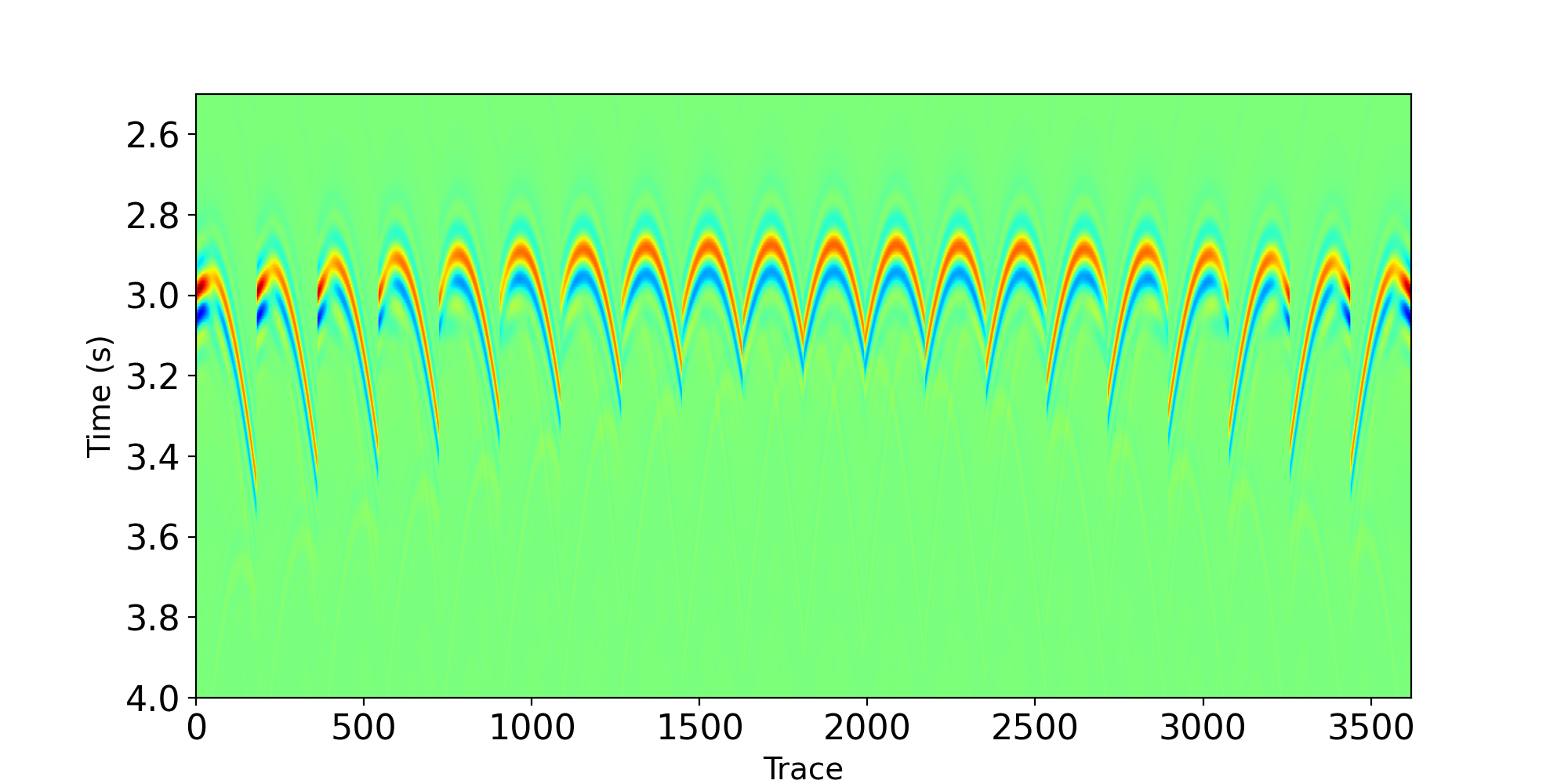}
\caption{Data generated using the Camembert
  model with the
 same source-receiver geometry as
  in Figure \ref{fig:dcwd20}: 20 sources located at $x=2000$ m, spaced
  150 m apart, starting at $z=500$ m; 181 receivers located at $x=6000$ m,
  spaced 20 m apart, starting at $z=200$ m.}
\label{fig:camcwd20}
\end{center}
\end{figure}

Figure \ref{fig:cam20uest0} shows the adaptive filter needed to match the data displayed in Figure \ref{fig:camcwd20} to the homogeneous medium data (Figure \ref{fig:chwd20}). The expected time lag of about 0.1 s is evident.


\begin{figure}[htbp]   
\begin{center}
\includegraphics[width=\textwidth]{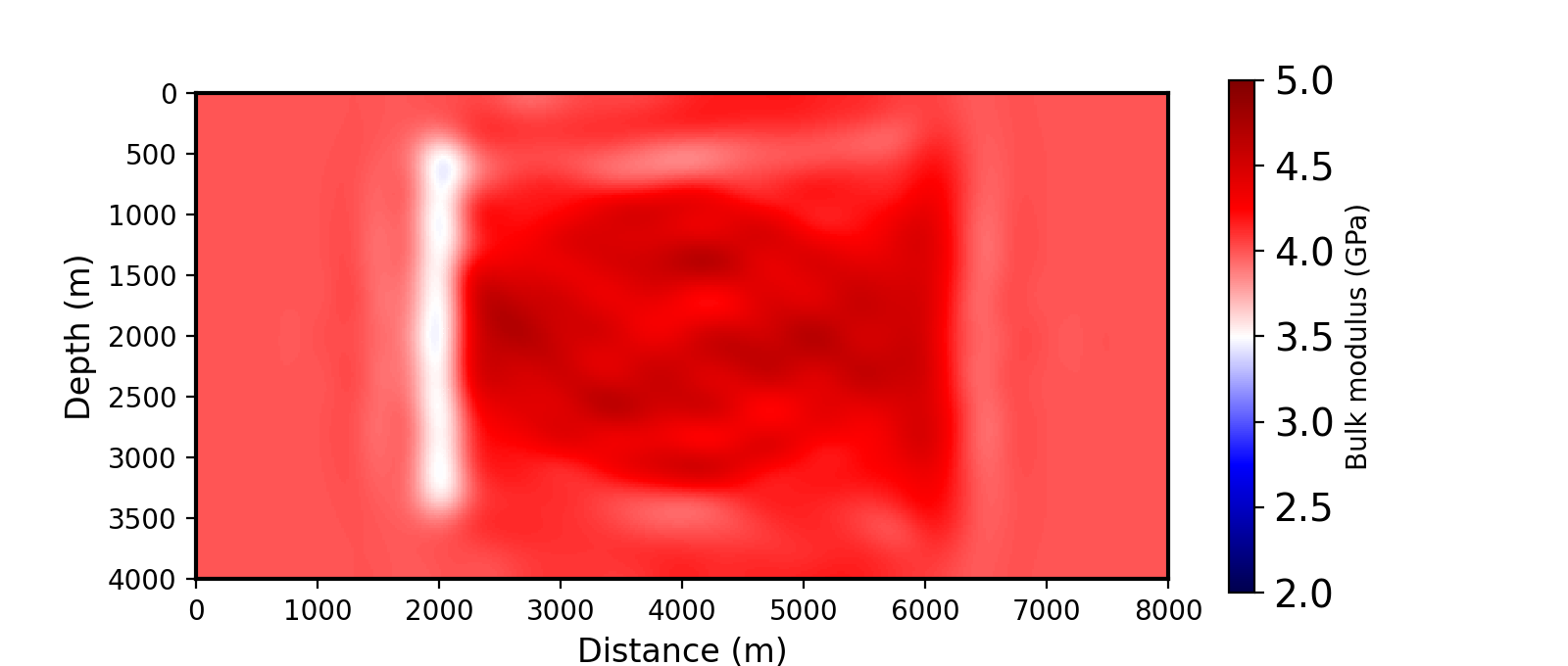}
\caption{Result of 12 MSWI
  iterations applied to the data in Figure \ref{fig:camcwd20} starting from
  the homogeneous initial model.}
\label{fig:cam20mestmswi}
\end{center}
\end{figure}

\begin{figure}[htbp]%
   \centering
   \subfloat[\centering ]{{\includegraphics[width=0.545\textwidth]{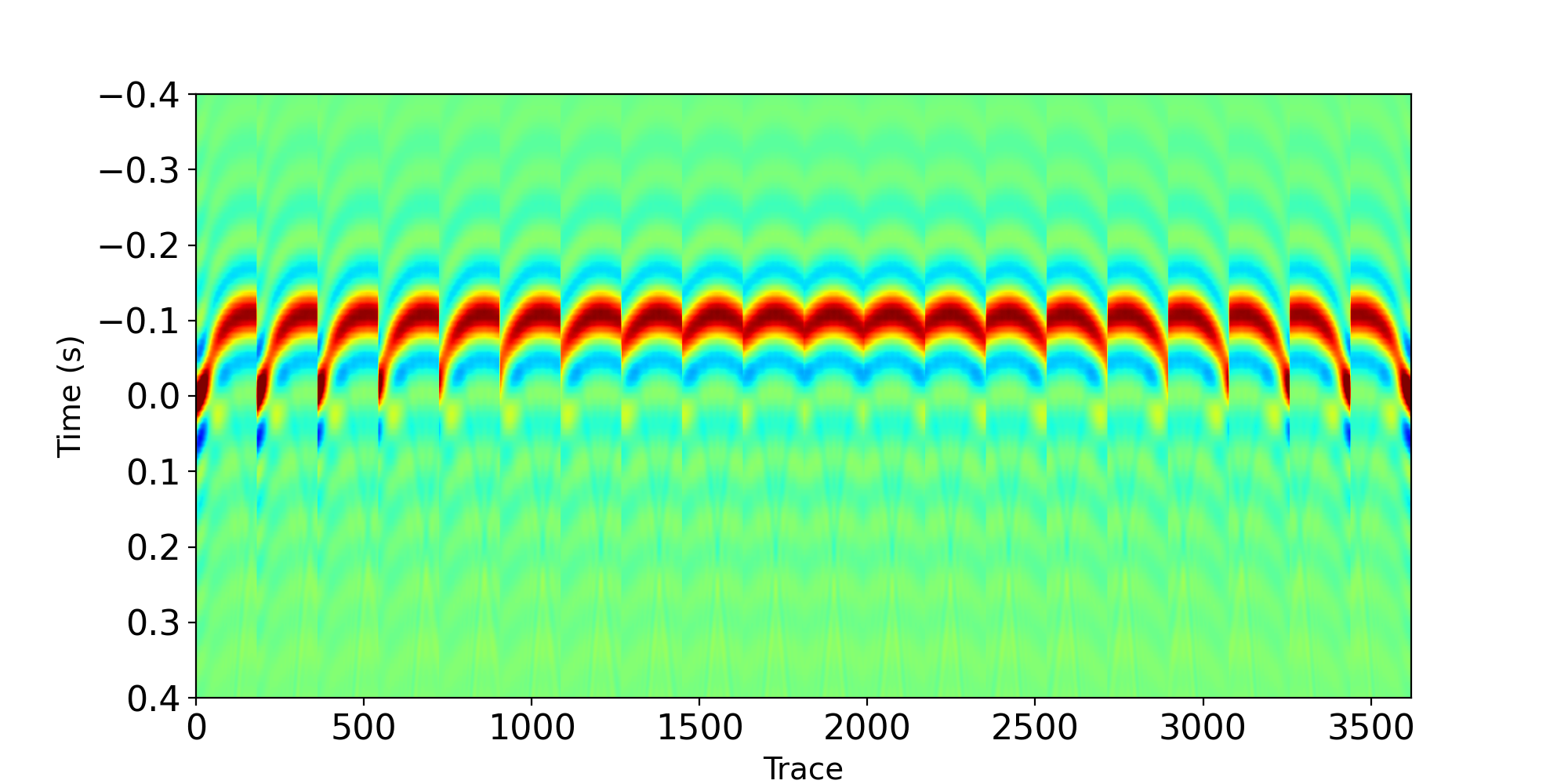}}\label{fig:cam20uest0}}%
   \hspace{-4em}
   \subfloat[\centering]{{\includegraphics[width=0.545\textwidth]{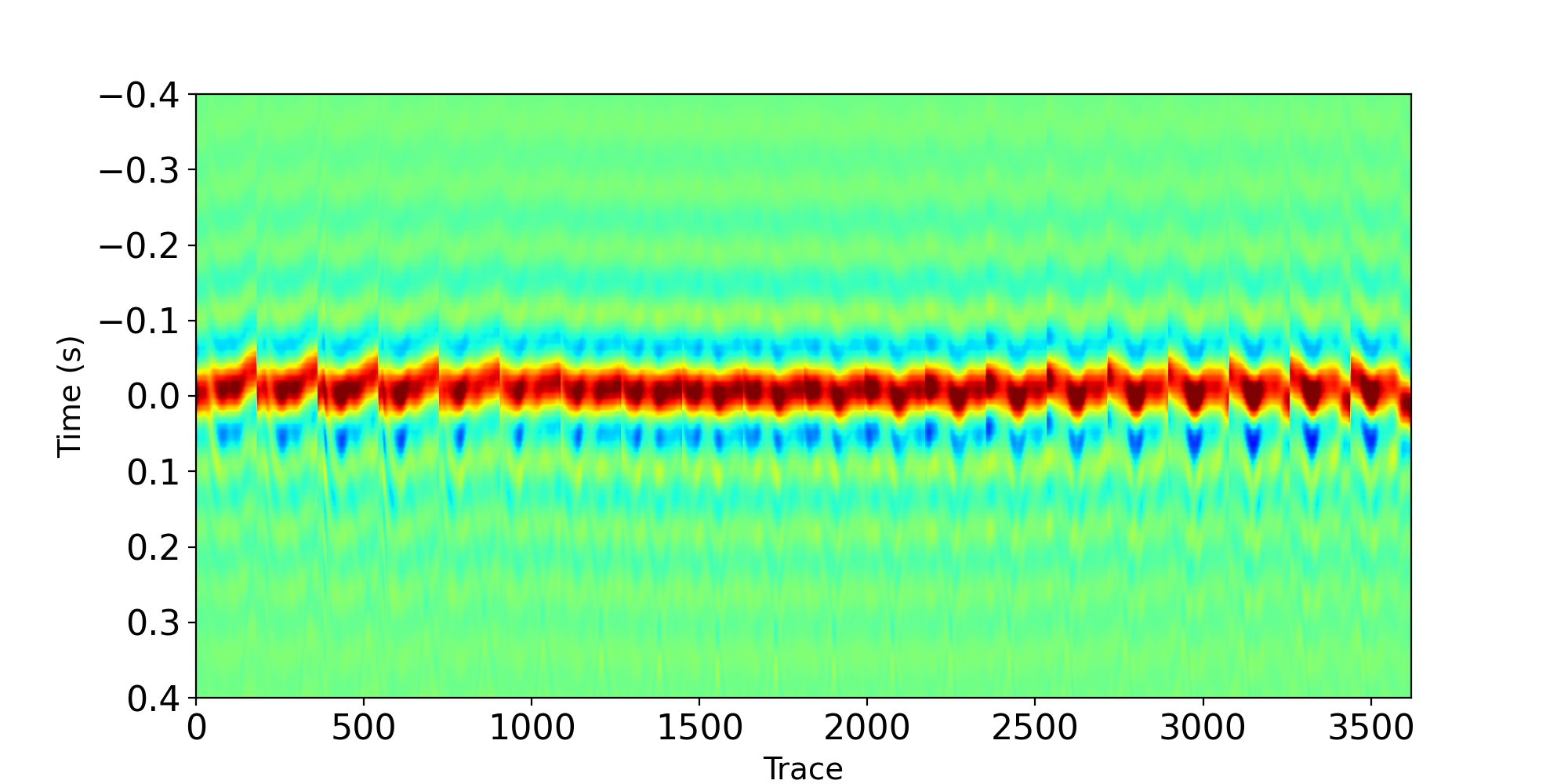} }\label{fig:cam20uestmswi}}%
   \caption{(a) Adaptive filter required to produce the data
of Figure \ref{fig:camcwd20} from corresponding homogeneous medium data (Figure
\ref{fig:chwd20}). (b) Adaptive filter estimate for MSWI
  inversion shown in Figure \ref{fig:cam20mestmswi}.}%
\end{figure}


As has been true for all of the examples shown in this section, iterative FWI starting at the homogeneous background model fails to generate a significant improvement in data fit with a modest number of iterations. On the other hand, 12 LBFGS iterations applied to the MSWI objective function, also starting at the homogeneous background model, yields the estimated bulk modulus field showns in Figure \ref{fig:cam20mestmswi}, which resembles a smoothed version of the target. It is smoothed because we use the same inverse weight operator in the bulk modulus model space (10 point smoothing in both directions, repeated once) as in the previous examples.  The adaptive filter returned at the final MSWI model (Figure \ref{fig:cam20uestmswi}) has greatly reduced the apparent time shift.

To accommodate the non-smooth nature of the target model in this example, we use an idea articulated most clearly by \cite{BarnierBiondi:23a}. Namely, for iterative inversion of wave data, the spatial frequency content of the trial models can be manipulated independently of the temporal frequency content of the data. The low spatial frequency components of trial models largely control the arrival times in predicted data, so it is natural to focus the inversion on them first. The previous examples have been based on smooth (low spatial frequency) target models, so determination of the low spatial frequency components was the whole story. Our version of the Camembert has high spatial frequency content. Rather than update the spatial frequency content during the extended inversion, as do \cite{BarnierBiondi:23a}, we observe that the extended (MSWI) inversion has produced a kinematically accurate model, and update spatial frequency content during FWI (as pointed out above, effectively the $\alpha \rightarrow \infty$ limit of MSWI). Figure 
\ref{fig:cam20mestmswifwifine} shows the result of 12 LBFGS iterations applied to the FWI objective, starting with the MSWI-estimated model in Figure \ref{fig:cam20mestmswi}, and with the same inverse weight operator as before, followed by 25 iterations of LBFGS with a ``hi-res" inverse weight operator (2 point moving average in both directions, repeated once). The final RMS error obtained is 
1.4 \% of the RMS error between the initial (homogeneous medium) predicted data (Figure \ref{fig:chwd20}) and the target data (Figure \ref{fig:camcwd20}). The apparent spatial resolution is much higher in the vertical direction than in the horizontal direction (an artifact of the source-receiver geometry). The data is much more sensitive to the change in travel time across the top and bottom edges of the Camembert than to those across the left and right edges.
\begin{figure}[htbp]   
\begin{center}
\includegraphics[width=\textwidth]{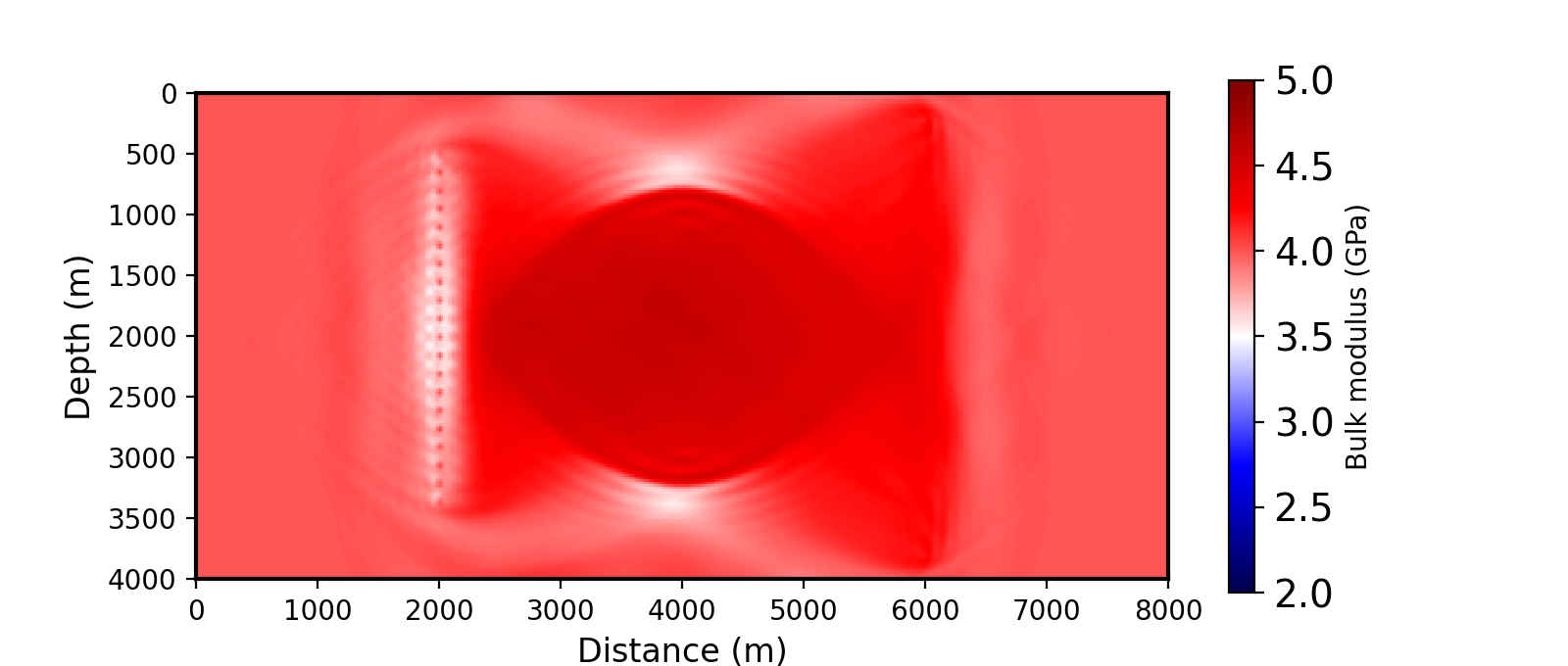}
\caption{``Hi-res'' bulk modulus produced by (1) 12 steps of LBFGS applied to MSWI objective, with (10 point moving average in both directions) inverse weight operator on bulk modulus followed by (2) 25 additional LBFGS steps, with (2 point moving average) inverse weight operator. 
  }
\label{fig:cam20mestmswifwifine}
\end{center}
\end{figure}



\subsection{Data with Additive Noise}


Seismic data contains substantial noise. Human activities such as road traffic and construction, and natural sources such as rivers, wind, and ocean waves can introduce noise in the data \cite[]{Shearer:2009}. Medical ultrasonics data often exhibits speckle noise 
\cite[]{KavandBekrani:24}. Finally, real world data are necessarily noisy 
because the underlying account of the physics is always simplified and hence incomplete. To test that our conclusions about MSWI are robust against the addition of modest amounts of noise, we present one example which includes 32\% (RMS) added noise. 

The 
example we describe in this subsection uses the target bulk modulus shown in Figure \ref{fig:m}. We generate noise by adding uniformly distributed random numbers (from $[-1,1]$ with mean zero and variance of $1/3$) to the homogeneous bulk modulus. The wave solution using this noisy model (and the same wavelet and buoyancy as in the circular lens experiment) is then added to the data shown in Figure \ref{fig:cwd20} to create the data for our experiments. 
In this example, we display single shot gathers 
in order to show clearly the effect of the additive noise.
The shot gather with source depth $z = 2000$ m appears in Figure \ref{fig:bncwlens20bncwd20}.
\begin{figure}[htbp]   
\begin{center}
\includegraphics[width=\textwidth]{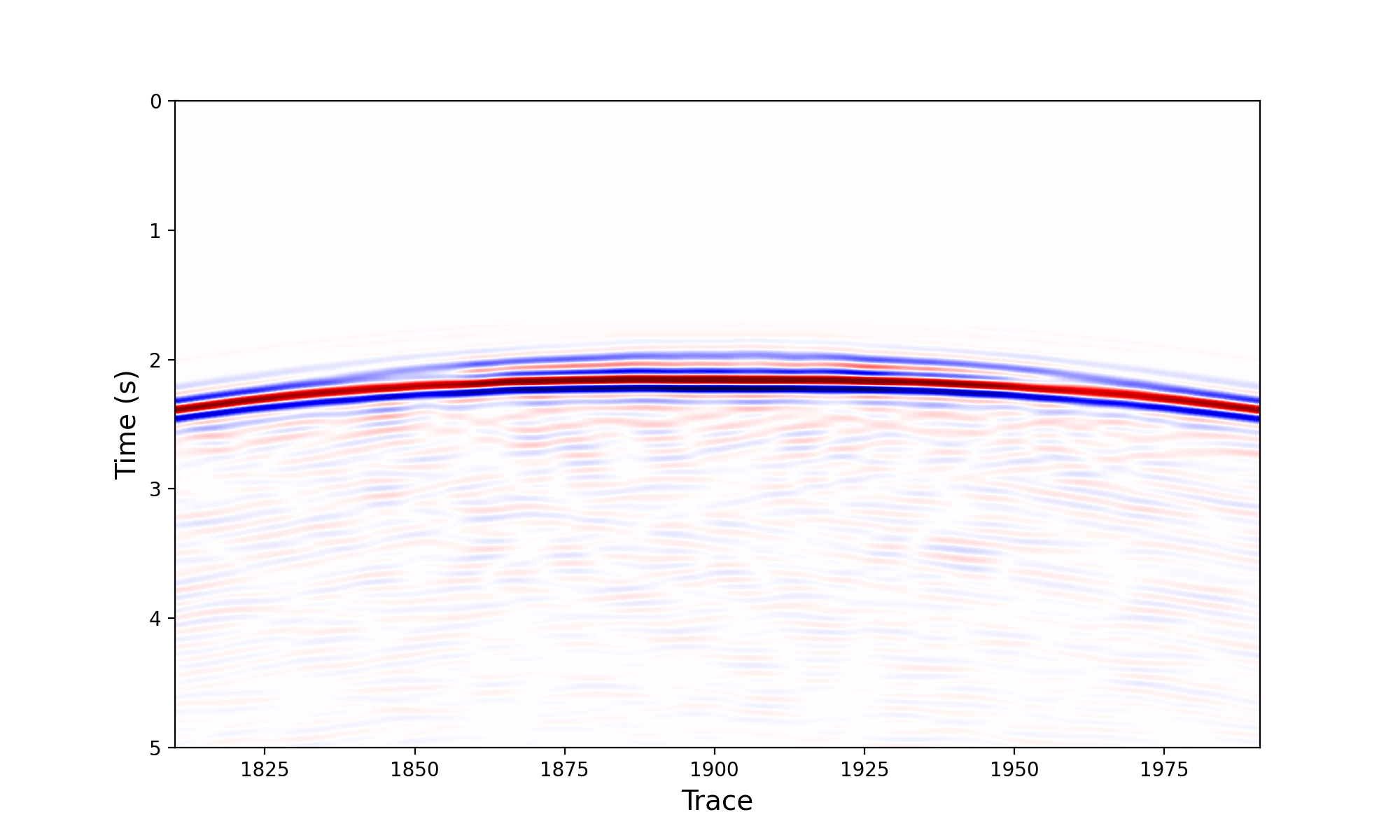}
\caption{Data for the circular lens model (Figure \ref{fig:m}) with approximately 32\% added noise. While there are 20 sources at $x=3000$ m,
here we show a representative set of data from one shot at depth $z$= 2000 m.}
\label{fig:bncwlens20bncwd20}
\end{center}
\end{figure}

We minimize the objective function $J_{\alpha,\sigma}[m_0, u; d]$ with respect to $u$ starting from data shown in Figure \ref{fig:chwd20}. The resulting adaptive filter $u_{\alpha,\sigma}[m_0;d]$ is shown in Figure \ref{fig:bncwlens20_covuest0}.
The initial value of the reduced MSWI objective
is $\approx 3.96 \times 10^{-2}$, and the initial gradient norm is $\approx 9.1 \times 10^{-5}$. After 12 iterations of LBFGS, the MSWI objective value has
decreased to $\approx 2.28 \times 10^{-3}$, and the norm of the gradient is
$\approx 9.1 \times 10^{-6}$. The resulting bulk modulus is shown in Figure
\ref{fig:bncwlens20_covmestmswi}, and 
the final adaptive filter is shown in Figure
\ref{fig:bncwlens20_covuestmswi}. 
Finally, starting from the MSWI estimated bulk modulus shown in Figure \ref{fig:bncwlens20_covmestmswi}, we performed 12 more iterations of LBFGS with the FWI objective function.
The FWI
objective value starts at $\approx 9.37$ and decreases to $\approx 2.87 \times
10^{-1}$. The resulting bulk modulus field is shown in Figure
\ref{fig:bncwlens20_covmestmswifwi}. 

\begin{figure}[htbp]%
   \centering
   \subfloat[\centering ]{{\includegraphics[width=0.5\textwidth]{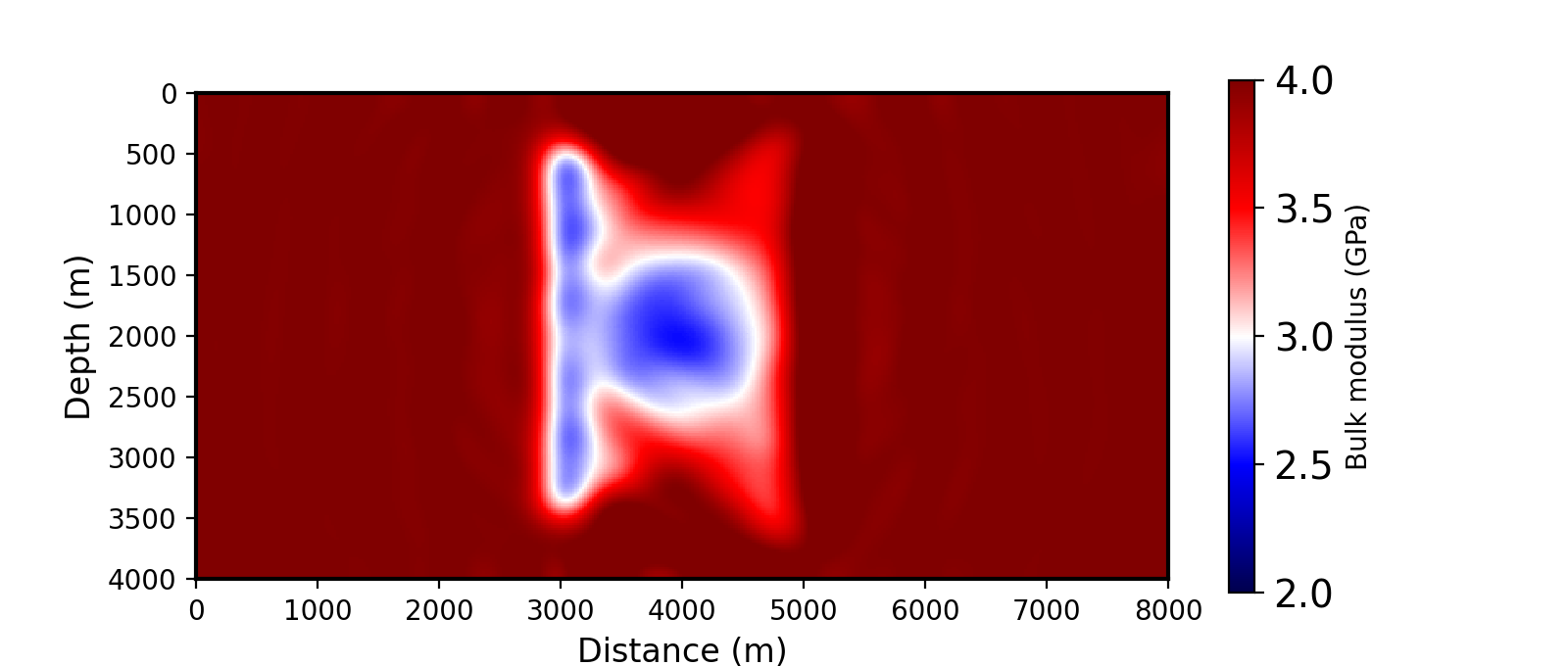}}\label{fig:bncwlens20_covmestmswi}}%
   \hspace{-2em}
   \subfloat[\centering]{{\includegraphics[width=0.5\textwidth]{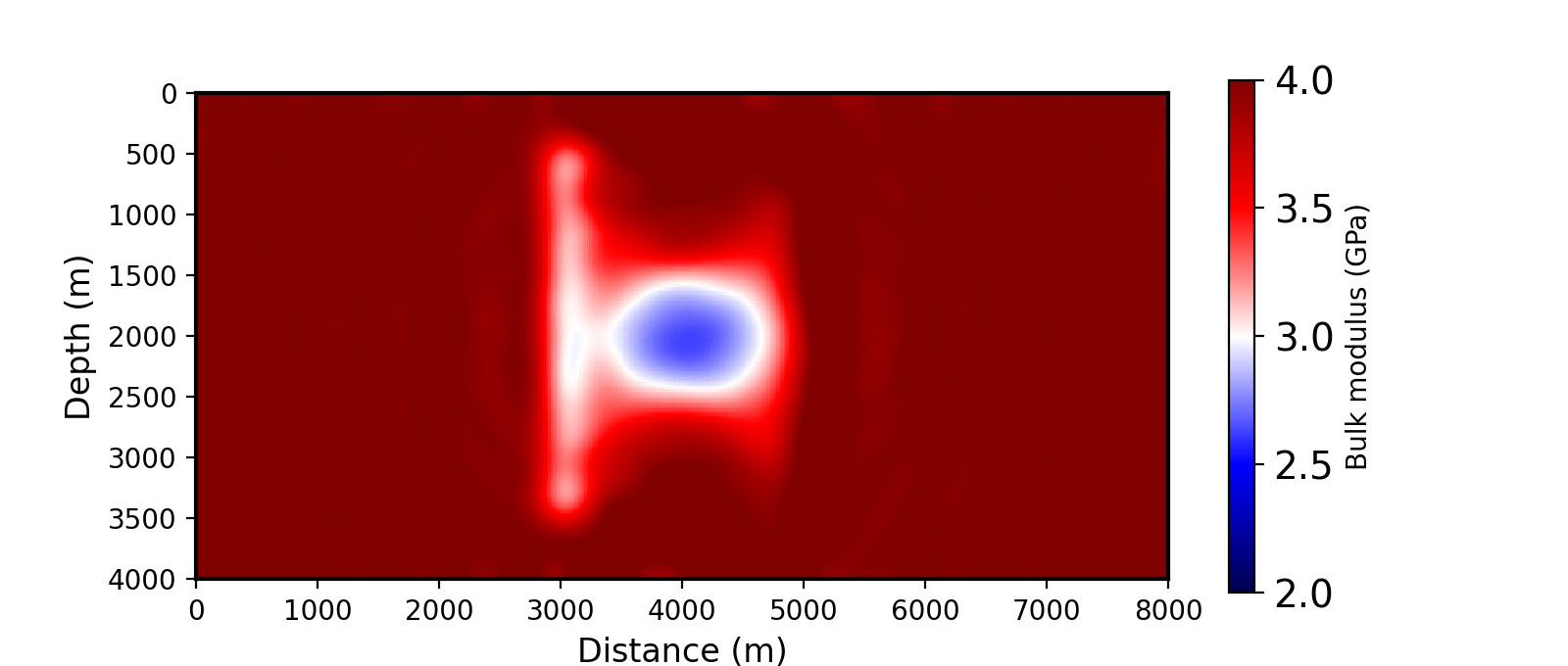} }\label{fig:bncwlens20_covmestmswifwi}}%
   \caption{(a) Bulk modulus produced by 12
  LBFGS steps to minimize the reduced MSWI objective 
  using the data (one shot of which is) shown in \ref{fig:bncwlens20bncwd20}. (b) Bulk modulus produced by 12
  LBFGS steps to minimize the FWI objective 
 starting at the MSWI
  result shown in Figure \ref{fig:bncwlens20_covmestmswi}.}%
\end{figure}

\begin{figure}[htbp]%
   \centering
   \subfloat[\centering ]{{\includegraphics[width=0.545\textwidth]{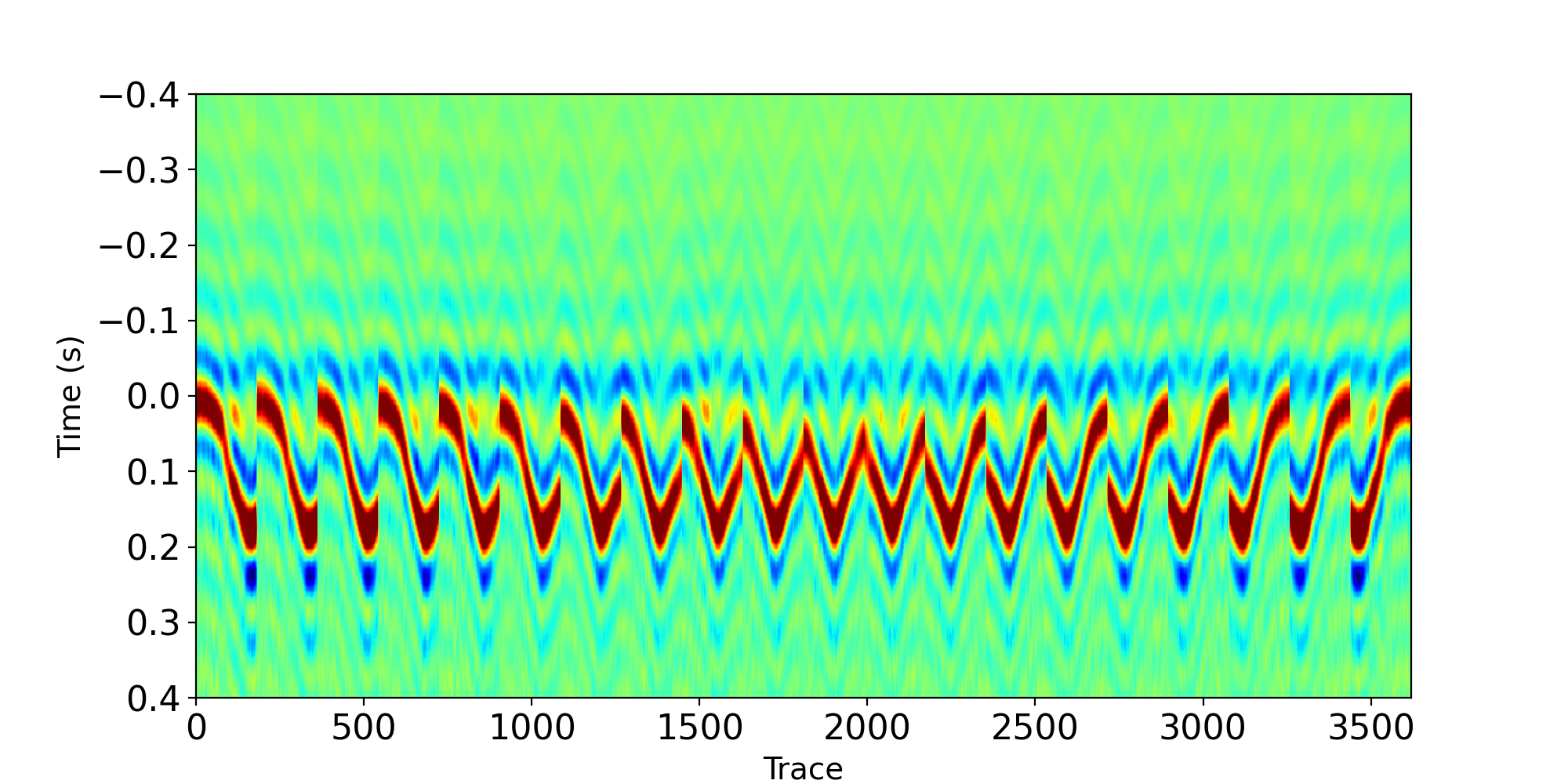}}\label{fig:bncwlens20_covuest0}}%
   \hspace{-4em}
   \subfloat[\centering]{{\includegraphics[width=0.545\textwidth]{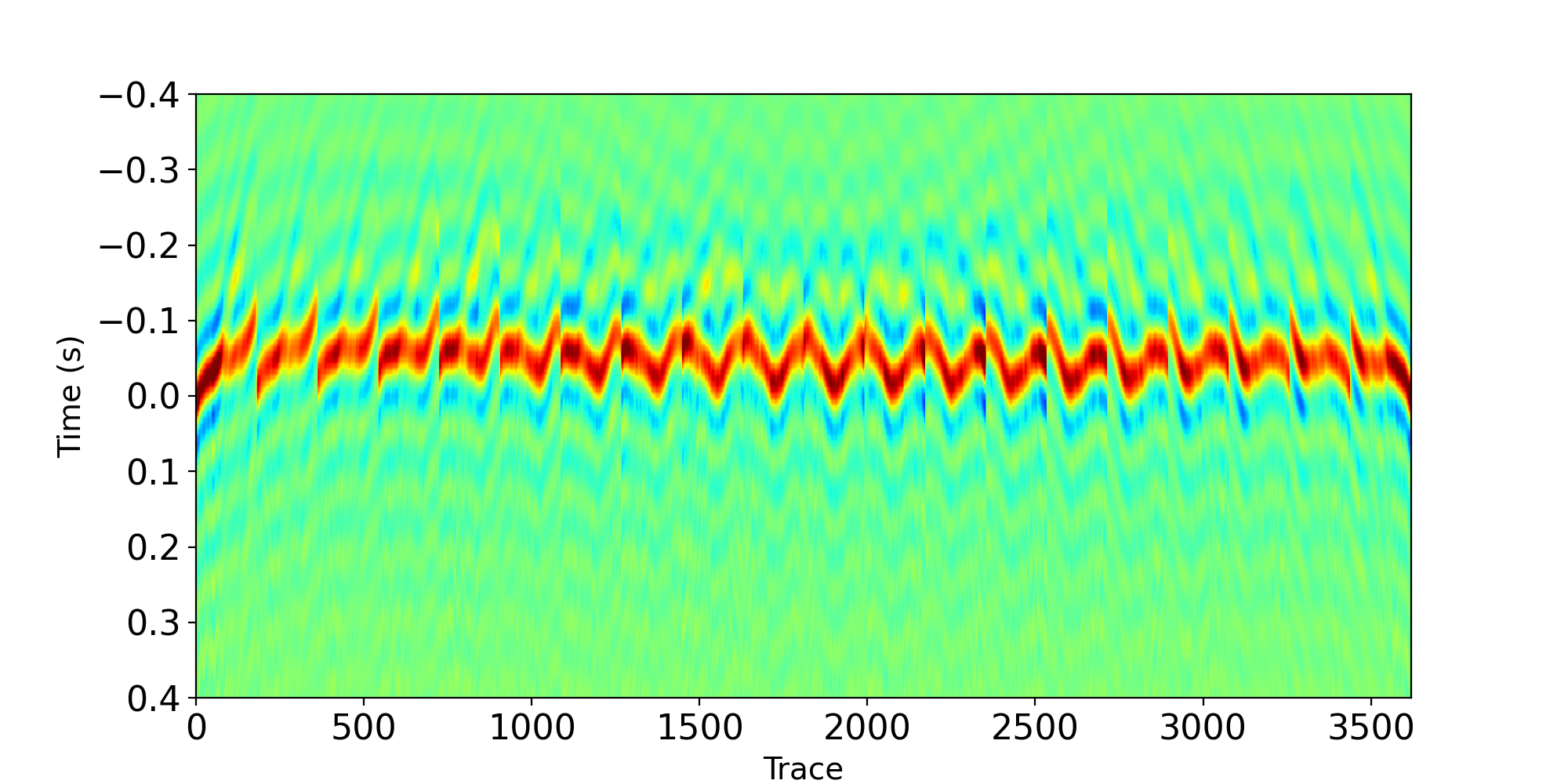} }\label{fig:bncwlens20_covuestmswi}}%
   \caption{(a) Adaptive filter to match
  initial model data (Figure \ref{fig:chwd20}) to circular lens
  data with 32\% added noise. (See Figure \ref{fig:bncwlens20bncwd20} for one shot of the noisy data.) Note considerable energy dispersion
  away from $t=0$. (b) Adaptive filter to match
  data for MSWI model
  to noisy circular lens
  data. 
  Note considerably less energy dispersion
  away from $t=0$ than is evident in Figure \ref{fig:bncwlens20_covuest0}.}%
\end{figure}

\section{Discussion}
In this section we will discuss some outstanding issues for the future including: the character of the FWI result, the choice of parameters in the
definition of the MSWI objective, and the application of MSWI to inverse problems
defined by other types of wave physics.

\subsection{FWI Failure Characterization} 
The FWI result shown in Figure \ref{fig:cwlens20mestfwi0} is certainly
unsatisfactory, as is evident from the residual plot (Figure
\ref{fig:cwlens20residmestfwi0}). 
We found that a modest increase (O(10)) in the number of LBFGS iterations did not significantly improve the data fit. We cannot conclude from this observation that the model shown in Figure \ref{fig:cwlens20mestfwi0} is near a stationary point or local minimizer of the FWI objective function or that sufficiently many more iterations would not approximate a data-fitting model. However,
even if orders of magnitude  
additional iterations {\em would} yield 
satisfactory data fit, such a vast amount of computational
work is simply not practical. Nor is it sensible, since alternatives
exist (FWI + MSWI) that accomplish the inversion goal with far less
work.
So, we don't know whether Figure \ref{fig:cwlens20mestfwi0} represents an
approximation to a non-global local minimizer or not, but it doesn't matter.

\subsection{Gradient Computation}
The reader will note that the theory cited above shows that the MSWI objective function is close to the travel time mean square error under some circumstances. It does not however assert that the {\em gradients} of these two types of objective function are close. This omission prevents any theoretical conclusion about a relation between stationary points of MSWI and travel time tomography objective functions. 

The compuation of the gradient via the Variable Projection gradient formula \ref{eqn:gradredfiltpen} also requires a better theoretical foundation. The transpose-derivative $DF^T$ of the modeling operator $F$ tends to enhance high-frequency components of model perturbations, relative to the output of the modeling operator itself. In the case presented here and similar approaches, the field to which $DF^T$ is applied is the result of iterative solution of a linear system, and the high frequency components of the iterates may not be assured to converge. In fact, the tangent approximation to this operator degrades in accuracy as maximum data frequency increases. The negative consequences of this phenomenon for accurate gradient computation, and possible fixes, have been discussed by several authors (see for example \cite{Symes:EAGE15,ChaurisFarshad:GEO23}). 

\subsection{Choosing Penalty Parameters}
The roles of parameters $\sigma$ and $\alpha$ are quite
different. Regularization, that is, avoidance of singularity in the
normal equation \ref{eqn:normal}, is the purpose of setting $\sigma >
0$. As formulated here, $\sigma$ is a dimensional parameter. It is
simple to make it dimensionless, and that should be done, so that a
``small'' choice is meaningful.

The penalty parameter $\alpha$ however plays a very different role. In
principle, $\alpha \rightarrow \infty$ should force the adaptive
filter to become ``physical'', which in this case means the identity
operator, i.e. with kernel $\delta(t)$ regardless of source and
receiver, so that the penalty problem is identical to FWI.
There are a variety of rules for setting penalty parameters like
$\alpha$. 
The discrepancy concept is one
of these. It is simple to implement, doesn't require a good initial value for problems like
MSWI, and steers the penalty problem towards its limit. 
The trade off is this approach requires the assertion of data error. The penalty parameter $\alpha$
is adjusted so that the data error term in \ref{eqn:redfiltpen} 
is close to the posited size (see \cite{Fu:Geo17b,SymesMinkoffChen:IP22}). 
One might note that we have traded one difficult-to-estimate 
parameter ($\alpha$) for another
(data error). However it is possible to turn this
construction around to estimate data error
\cite[]{ChenSymesMinkoff:IMAGE22}. In any case, data error has an
obvious meaning, whether it is a known quantity or not, whereas
$\alpha$ is merely a Lagrange multiplier.

\subsection{Complex Wavefronts}
As has been emphasized several times, wavefront complexity can prevent
successful application of inversion methods based on source-receiver
extension. The theoretical link between MSWI and similar methods, on the one hand, and travel time tomography, on the other, is broken in the presence of multiple ray paths between source and receiver, and numerical examples conform to the theory as already shown in the second oblate lens example (see also \cite{Symes:94c,HuangSymes:Geo17,Symes:24a}). 
However, the
source-receiver extension is not the only possible route to FWI
modification by artificially expanding the definition of energy
source. \cite{HuangNammourSymesDollizal:SEG19} overview modifications
of FWI based on various source extensions; some more recent advances
are described by
\cite{MetivierBrossier:SEG20,PladysBrossierLiMetivier:GEO21,LiAlkhalifah:21,Yongetal:GJI23,Opertoetal:GEO23,Symes:23}.
Numerical examples suggest that some of these extensions may avoid
cycle-skipping for transmitted wave data with complex wavefronts.
For many of these
approaches, the missing piece is a theoretical justification like that
summarized in equations \ref{eqn:tomo1}-\ref{eqn:tomo3} for MSWI.

\subsection{Applying MSWI to Other Types of Physics and Data}
Acoustics is seldom exactly the right choice of physics for mechanical
vibration modeling. In fact in many settings some variant of elastodynamics,
perhaps coupled with acoustics in some regions, is the right
choice. Elastic waves travel at several speeds, intrinsically
providing multiple wavefronts. Extension of MSWI 
to elastodynamic inverse problems, either by data mode
separation or some more intrinsic construction, would appear to be a
very important next step, as would application to electromagnetic imaging.

Sensors (and energy sources) are seldom punctual and/or
isotropic. Actual sources encountered in lab and field experiments can exhibit
non-trivial radiation patterns or even a physical extent which is non-negligible
on the wavelength scale. Inclusion of more realistic source modeling
also seems a very important next step
for testing the limits of source-receiver extension.

Finally, the assumption of smooth material parameter variation is
unlikely to be justified in actual physical media. Sedimentary rocks,
human tissue, and manufactured objects all exhibit juxtaposition of
structure at all scales, from long scales modeled by smooth variation,
to sub-wavelength oscillations. Most likely in settings where
ballistic transmitted waves dominate the data energetics, techniques like MSWI can be applied
to
gain a low-resolution image of model parameters like wave velocity or bulk modulus, etc. However
theoretical results in this direction are absent, as are results on
the related question of applicability of algorithms like MSWI to
predominantly reflected wave data.

\section{Conclusions}
Matched Source Waveform Inversion is one of several extension methods that can avoid the cycle skipping phenomenon afflicting standard least-squares (full waveform) inversion in certain special cases.  The adaptive filter in MSWI helps align the predicted and observed data, which allows for larger model adjustments than can be made using local-gradient based optimization on the FWI objective. In this paper we provide a catalogue of numerical experiments that illustrate several features of MSWI. All of these experiments use relatively simple (and similar) heterogeneous velocity models with a single anomaly embedded within a homogeneous background medium. They all assume
a single vertical line of sources on the left of the anomaly and another vertical line of receivers on the right. Despite their simplicity, these velocity models illustrate how and when MSWI can provide a good starting model for FWI. In a related theory paper, MSWI applied to single-arrival data is shown to be equivalent to travel-time tomography which does not suffer from cycle skipping. We first show an experiment with single arrival data that 
illustrates failure of FWI to converge to a reasonable velocity when the inversion is started from a homogeneous initial model. MSWI, on the other hand, succeeds for this same experiment. Further, this MSWI result can be used as an initial guess for FWI allowing further refinement of the model. However, when
the model and acquisition geometry of the experiment change so that the data includes multiple arrivals, our numerical examples indicate that MSWI is no
longer able to find a reasonable velocity from a homogeneous initial model. Additionally we show that with moderate levels of added
noise, MSWI is still able to converge to a reasonable velocity model from a homogeneous starting model. 

\section{Acknowledgements}
The authors gratefully acknowledge support from the sponsors of the UT Dallas {\it 3D+4D Seismic FWI} research consortium.

\bibliographystyle{seg}
\bibliography{masterref}
\end{document}